\documentclass[journal]{IEEEtran}
\usepackage[utf8]{inputenc}
\usepackage{color}
\usepackage{array}
\usepackage{verbatim}
\usepackage{float}
\usepackage{amsmath}
\usepackage{amsthm}
\usepackage{amssymb}
\usepackage{graphicx}
\usepackage{longtable}
\usepackage[unicode=true,
 bookmarks=false,
 breaklinks=false,pdfborder={0 0 1},colorlinks=false]
 {hyperref}
\hypersetup{
 colorlinks,bookmarksopen,bookmarksnumbered,citecolor=blue,urlcolor=blue}
\usepackage[]{apacite}

\usepackage{lipsum}
\usepackage{mathtools}
\usepackage{cuted}

\floatstyle{ruled}
\newfloat{algorithm}{tbp}{loa}
\providecommand{\algorithmname}{Algorithm}
\floatname{algorithm}{\protect\algorithmname}

\makeatletter
\let\oldforeign@language\foreign@language
\DeclareRobustCommand{\foreign@language}[1]{%
	\lowercase{\oldforeign@language{#1}}}

\let\oldforeign@language\foreign@language
\DeclareRobustCommand{\foreign@language}[1]{%
	\lowercase{\oldforeign@language{#1}}}

\ifCLASSINFOpdf
\else
\fi

\newcommand{\MYfooter}{\smash{
		\hfil\parbox[t][\height][t]{\textwidth}{\centering
			\thepage}\hfil\hbox{}}}

\def\ps@IEEEtitlepagestyle{%
	\def\@oddhead{\parbox[t][\height][t]{\textwidth}{\centering \scriptsize
			Personal use of this material is permitted. Permission from the author(s) and/or copyright holder(s), must be obtained for all other uses. Please contact us and provide details if you believe this document breaches copyrights.\\
			\noindent\makebox[\linewidth]{}
		}\hfil\hbox{}}%
	\def\@evenhead{\scriptsize\thepage \hfil \leftmark\mbox{}}%
	\def\@oddfoot{\parbox[t][\height][l]{\textwidth}{
			\vspace{-20pt}{\rule{\textwidth}{0.4pt}}\\ \footnotesize\underline{To cite this article:}
			{\bf{\footnotesize\textcolor{red}{H. A. Hashim, S. El-Ferik, and F. L. Lewis, ”Neuro-adaptive Cooperative Tracking Control with Prescribed Performance of Unknown Higher-order Nonlinear Multi-agent Systems,” International Journal of Control, vol. 92, no. 2, pp. 445-460, 2019.}}} doi: \href{https://doi.org/10.1080/00207179.2017.1359422}{10.1080/00207179.2017.1359422}\\
			\noindent\makebox[\linewidth]
		}\hfil\hbox{}}%
	\def\@evenfoot{\MYfooter}}

\makeatother
\newtheorem{defn}{Definition}

\newtheorem{lem}{Lemma}

\newtheorem{thm}{Theorem}
\newtheorem{rem}{Remark}

\newtheorem{assum}{Assumption}

\begin{document}
%\onecolumn
%\noindent\rule{18.1cm}{2pt}\\
%\underline{To cite this article:}
%{\bf{\textcolor{red}{H. A. Hashim, S. El-Ferik, and F. L. Lewis, "Neuro-adaptive Cooperative Tracking Control with Prescribed Performance of Unknown Higher-order Nonlinear Multi-agent Systems," International Journal of Control, vol. 92, no. 2, pp. 445-460, 2019.}}}\\
%\noindent\rule{18.1cm}{2pt}\\
%
%\noindent{\bf The published version (DOI) can be found at: \href{http://dx.doi.org/10.1080/00207179.2017.1359422}{10.1080/00207179.2017.1359422} }\\
%
%\vspace{40pt}\noindent Please note that where the full-text provided is the Author Accepted Manuscript or Post-Print version this may differ from the final Published version. { \bf To cite this publication, please use the final published version.}\\
%
%
%\textbf{
%	\begin{center}
%		Personal use of this material is permitted. Permission from the author(s) and/or copyright holder(s), must be obtained for all other uses, in any current or future media, including reprinting or republishing this material for advertising or promotional purposes.\vspace{60pt}\\
%	\end{center}
%\vspace{360pt}
%%	\begin{flushleft}
%%		\underline{Publication information:}\\
%% Date of Submission: March 2017.\\
%%Date of Acceptance: July 2017.\\
%%Available Online: August 2017.
%%\end{flushleft} 
%}
%\footnotesize{ \bf
%	\vspace{20pt}\noindent Please contact us and provide details if you believe this document breaches copyrights. We will remove access to the work immediately and investigate your claim.
%} 
%
%\normalsize
%
%
%\twocolumn
\title{Neuro-adaptive Cooperative Tracking Control with Prescribed Performance of Unknown Higher-order Nonlinear Multi-agent Systems}

\author{Hashim~A.~Hashim$^*$, Sami~El-Ferik, and~Frank L. Lewis% <-this % stops a space
\thanks{$^*$Corresponding author, H. A. Hashim is with the Department of Electrical and Computer Engineering,
University of Western Ontario, London, ON, Canada, N6A-5B9, e-mail: hmoham33@uwo.ca}% <-this % stops a space
\thanks{S.~El-Ferik is with the Department of Systems Engineering, King Fahd University of Petroleum and Minerals, Dhahran, 31261, Saudi Arabia.}% <-this % stops a space
\thanks{F. L. Lewis is with the Department of Electrical and Computer Engineering, UTA Research Institute, The University of Texas at Arlington 7300 Jack Newell Blvd. S, Ft. Worth, Texas 76118.}% <-this % stops a space
}

%\author{Hashim~A.~Hashim$^*$,~\IEEEmembership{~Member, IEEE}, Lyndon J. Brown, and~Kenneth McIsaac,~\IEEEmembership{~Fellow, IEEE}% <-this % stops a space
%\thanks{$^*$Corresponding author, H. A. Hashim, L. J. Brown and K. McIsaac are with the Department of Electrical and Computer Engineering,
%University of Western Ontario, London, ON, Canada, N6A-5B9, e-mail: hmoham33@uwo.ca, lbrown@uwo.ca and kmcisaac@uwo.ca.}}

\markboth{--,~Vol.~-, No.~-, \today}{Hashim \MakeLowercase{\textit{et al.}}: Nonlinear Stochastic Attitude Filter on the Special Orthogonal Group}
\markboth{}{Hashim \MakeLowercase{\textit{et al.}}: Nonlinear Stochastic Attitude Filter on the Special Orthogonal Group}

\maketitle

\begin{abstract}
 This paper is concerned with the design of a distributed cooperative synchronization controller for a class of
higher-order nonlinear multi-agent systems. The objective is to achieve synchronization and satisfy a predefined time-based performance. Dynamics of the agents (also called the nodes) are assumed to be unknown to the controller and are estimated using Neural Networks. The proposed robust neuro-adaptive controller drives different states of nodes systematically to synchronize with the state of the leader node within the constraints of the prescribed performance. The nodes are connected through a weighted directed graph with a time-invariant topology. Only few nodes have access to the leader. Lyapunov-based stability proofs demonstrate that the multi-agent system is uniformly ultimately bounded stable. Highly nonlinear heterogeneous networked systems with uncertain parameters and external disturbances were used to validate the robustness and performance of the new novel approach. Simulation results considered two different examples: single-input single-output and multi-input multi-output, which demonstrate the effectiveness of the proposed controller.
\end{abstract}


\IEEEpeerreviewmaketitle{}

\section{Introduction}

\IEEEPARstart{T}{he} use of collaborative autonomous robotic vehicles allows for greater flexibility and capacity as well as higher performance in areas such as surveillance, inspection, space explorations, communication, sensor deployment and many others. Multi-agent systems (MAS) distribute work in a logical manner and exchange information via self-formed local network and, hence, they are often called nodes. The network is named a communication graph formed by a set of nodes and the communication lines between different nodes are called edges. The graph can be directed or undirected. An undirected graph allows the information to flow in both directions. The connected nodes of such a graph own similar characteristics. On a directed graph or a digraph the direction of the information flow is fixed. The direction is pointed from one node to another indicating how the information flows from one node to its neighbors. Moreover, the structure of the network can be fixed or variable.\\
The control of such multi-agent systems faces several practical as well as theoretical challenges (see for instance \cite{olfati-saber_consensus_2004}).  In particular, dynamics of the node can be nonlinear and unknown, the network bandwidth capacity is limited and may suffer from variable delays and loss of packets, the operating environment is changing and complex with presence of noise, the embedded computational resources are limited, etc. In the literature, several studies addressed either cooperative regulation problem, called consensus, or cooperative tracking problem, known as synchronization (see for example \cite{fax_information_2004}). Recently, several control methods for high-order non-linear multi-agent systems have been proposed. Synchronization of passive nonlinear systems has been considered in \cite{chopra_passivity-based_2006} while distributive tracking problem of node consensus has been studied extensively such as \cite{lewis_cooperative_2013}, \cite{olfati-saber_consensus_2007}, \cite{liao2016cooperative} and \cite{zhang_adaptive_2012}. Work of \cite{das_distributed_2010} and \cite{cao_distributed_2012} studied cooperative tracking control for single node representing a single-input single-output (SISO) system with high order dynamics. Due to unknown dynamics, \cite{das_distributed_2010} and \cite{zhang_adaptive_2012} proposed a neuro-adaptive distributed control for heterogeneous agents connected through a digraph. \cite{das_distributed_2010} considered single integrator agents and later on \cite{zhang_adaptive_2012}, high order affine systems described in Brunovsky form and connected through a directed graph have been addressed. The authors assumed that the input function $g_{i}\left(\cdot\right)$ is equal to one for each agent $i=1, \ldots, N$.\\
In all previous studies the input function was assumed to be known. On the other hand, adaptive distributed tracking control of affine systems has been studied assuming unknown input function by \cite{theodoridis_direct_2012} and extended in \cite{el-ferik_neuro-adaptive_2014}. Also, consensus with Saturation and Dead-zone was examined in \cite{shen2016adaptive, shen2016output}. \cite{theodoridis_direct_2012} approximated the unknown nonlinear dynamics and input functions using a neuro-adaptive fuzzy and defined the output membership functions by a set of offline trials. All these previous studies mainly focused on ultimate stability of the error response. Most of the proposed controllers for highly nonlinear systems guarantee that the consensus tracking error is upper bounded due to uncertainties in dynamics and external disturbances. 
Consensus in error has been proven to be ultimately uniformly bounded and to converge into a residual set having a size that depends on some unknown  but bounded sets. However, bounded sets represent uncertainties in dynamics and external disturbances.
Therefore, it is almost impossible to make the prediction of transient performance as well as steady-state behavior analytically  \cite{bechlioulis_robust_2008}.\\
On the other hand, designing a cooperative adaptive control for a group of agents satisfying prescribed performance function (PPF) has some advantages. PPF forces the output error to begin within large set and  steer systematically into an arbitrarily small set satisfying a known measure \cite{bechlioulis_robust_2008} and \cite{Hashim2017adaptive}. Under prescribed performance, the error should display some measure of dynamic features. For instance, convergence rate should obey a predefined value and a maximum value of overshoot or undershoot is not exceeding a given range. In addition to having the error dynamically bounded, prescribed performance-based controller for cooperative adaptive control is capable of reducing the control effort and improving its robustness. Upper and lower bounds of PPF should be defined appropriately in order to provide smooth tracking error with prescribed convergence. Neuro adaptive control with PPF for strict feedback linearizable systems has been presented by \cite{bechlioulis_robust_2008}. Since then, several papers developed neuro adaptive control with prescribed performance approximating the unknown nonlinearities and disturbances through linearly parametrized neural network (see for instance \cite{bu2016robust}, \cite{yang_adaptive_2015} and \cite{el2017neuro}). A model reference adaptive control with PPF has been proposed to avoid defining neural weights via trial and error methods \cite{mohamed_improved_2014}. Most of the studies considered only single autonomous systems. However, just recently \cite{Hashim2017adaptive} considered networked graph and proposed an adaptive cooperative control with prescribed performance for a first order node dynamics with unknown nonlinearities. \cite{zhang2016distributed} addressed the problem of distributed output feedback consensus tracking control for leader following nonlinear multi-agent systems in strict-feedback form  with PPF requirement. A similar work has been proposed by \cite{shahvali2016cooperative}.\\
%The goal of this paper is to extend the results of \cite{zhang_adaptive_2012} with a control design based on a decentralized neuro-adaptive cooperative control for a fleet of autonomous systems linked through a communication digraph and required to perform within a predefined set of dynamical requirements. Like \cite{zhang_adaptive_2012}, the system considered is in Brunovsky form. However, the different synchronization errors between all these nodes are expected to have dynamic responses within some performance constraints determined by the designer. the dynamics of each node is subject to model uncertainties and to unknown but bounded external disturbances. The paper can also be seen as an extension to \cite{shahvali2016cooperative}, where the author rightly identified the need to address the control of highly nonlinear systems with prescribed performance. 
Indeed, the present proposed control scheme is developed using prescribed performance to satisfy transient and steady-state dynamic performance for each node's state through synchronization error. The data exchange between nodes is carried out according to a given directed graph. Neural Network is used to estimate the unknown nonlinear dynamics. In addition, this paper considers the original prescribed performance scheme presented by \cite{bechlioulis_robust_2008}. Hence, the interactions between all nodes are considered in the consensus algorithm to track the leader trajectory and guarantee stable non-oscillatory dynamics.\\
The rest of the paper is organized as follows. Section \ref{Sec2} presents graph theory preliminaries and math notations. In Section \ref{Sec3}, problem formulation, the associated local error synchronizations and prescribed performance characteristics are formulated. Section \ref{Sec4} develops the control law in order to prove stability of the directed connected graph and satisfy prescribed performance characteristics. Section \ref{Sec4} also presents the neural approximation and stability of the control design of distributed agents based on neural approximation. Section \ref{Sec5} illustrates results which guarantee effectiveness and robustness of the proposed control for SISO and MIMO problems. Finally, conclusion and future directions of research are given in Section \ref{Sec6}.\\

\section{Preliminaries}\label{Sec2}
\subsection{Mathematical Identities}
Throughout this paper, the set of real numbers is denoted as $\mathbb{R}$;
$n$-dimensional vector space as $\mathbb{R}^{n}$; the space span
by $n\times m$ matrix as $\mathbb{R}^{n\times m}$; identity matrix
of order $m$ as $\mathbb{I}_{m}$; absolute value as $|\cdot|$.
For $x\in\mathbb{R}^{n}$, the Euclidean norm is given as $\left\Vert x\right\Vert =\sqrt{x^{\top}x}$
and matrix Frobenius norm is given as $\left\Vert \cdot\right\Vert _{F}$. For any $x_{i}\in\mathbb{R}^{n}$ we have $x_{i} = \left[x_{i}^{1},\ldots,x_{i}^{n}\right]^{\top}$ for $i=1,\ldots,N$ and for $x^{j}\in\mathbb{R}^{N}$ we have $x^{j} = \left[x_{1}^{j},\ldots,x_{N}^{j}\right]$ for $j=1,\ldots,n$.
Trace of associated matrix is denoted as ${\rm Tr}\left\{ \cdot\right\} $, ${\rm diag}\left\{\cdot\right\}$ denotes the diagonal of associated matrix,
$\mathcal{N}$ is the set $\{1,...,N\}$, and ${\bf \underline{1}}_{N}$
is a unity vector $[1,\ldots,1]^{\top}\in\mathbb{R}^{N}$. $A$ is
said to be positive definite if $A>0$ for $A\in\mathbb{R}^{n\times n}$;
$A\geq0$ indicates positive semi-definite; $\sigma\left(\cdot\right)$
is the set of singular values of a matrix with maximum value $\bar{\sigma}\left(\cdot\right)$
and minimum value $\underline{\sigma}\left(\cdot\right)$. Finally, $\otimes$ denotes the Kronecker product.\\
\subsection{Basic graph theory}
A graph is denoted by $\mathcal{G} = \left(\mathcal{V}, \mathcal{E}\right)$ with a nonempty finite set of nodes (or vertices)  $\mathcal{V} = \left\{\mathcal{V}_1, \mathcal{V}_2, \ldots, \mathcal{V}_n\right\}$, and a set of edges (or arcs) $\mathcal{E}\subseteq  \mathcal{V}\times \mathcal{V}$. $\left(\mathcal{V}_{i},\mathcal{V}_j\right) \in \mathcal{E}$ if there is an edge from node $i$ to $j$. 
Topology of a weighted graph is described by the adjacency matrix $A=\left[a_{i,j}\right]\in \mathbb{R}^{N\times N}$  with weights $a_{i,j} > 0$ if $\left(\mathcal{V}_{j}, \mathcal{V}_{i}\right) \in \mathcal{E}$: otherwise $a_{i,j} = 0$. 
Throughout the paper, a directed graph is called diagraph. 
Also, the topology is fixed where $A$ is time-invariant and the self-connectivity element $a_{i,i} = 0$. A graph can be directed or undirected.
The weighted in-degree of a node $i$ is given by the sum of $i$-th row of $A$, i.e., $d_{i}=\sum_{j=1}^{N} a_{i,j}$. Also, the diagonal in-degree matrix is $D ={\rm  diag}\left(d_1,\ldots, d_N\right)\in \mathbb{R}^{N\times N}$ and the graph Laplacian matrix $L = D-A$.
The set of neighbors of a node $i$ is $N_{i} = \left\{j|\left(\mathcal{V}_j\times \mathcal{V}_{i}\right)\in \mathcal{E}\right\}$. If node $j$ is a neighbor of node $i$, then node $i$ can get information from node $j$ , but not necessarily vice versa. For undirected graph, neighborhood is a mutual relation. A direct path from node $i$ to node $j$ is a sequence of successive edges in the form $\left\{\left(\mathcal{V}_{i},\mathcal{V}_k\right),\left(\mathcal{V}_k,\mathcal{V}_l\right),\ldots, \left(\mathcal{V}_m,\mathcal{V}_j\right)\right\}$. 
If there is a node such that there is a directed path from one node to every other node in the graph, then the diagraph has a spanning tree. 
If for any ordered pair of nodes $\left[\mathcal{V}_{i},\mathcal{V}_j\right]$ with $i\neq j$, then a diagraph is strongly connected and there is a directed path from node $i$ to $j$ \cite{ren_distributed_2008}.

\section{Problem Formulation in Prescribed Performance} \label{Sec3}
Let the nonlinear dynamics of the $i$th node be given by
\begin{equation}
\label{eq:equa1}
\begin{aligned}
\dot{x}_{i}^{1} & = x_{i}^{2}\\
\dot{x}_{i}^{2} & = x_{i}^{3}\\
& \vdots\\
\dot{x}_{i}^{M_{p}} & = f_{i}\left(x_{i}\right) + G_{i}u_{i}\\
y_{i} & = x_{i}^{1}
\end{aligned}
\end{equation}
where $x_{i}^{m_{p}} \in \mathbb{R}^{P}$ is the $m_{p}$th-state node of $i$ where $x_{i} = \left[x_{i}^{1}, \ldots, x_{i}^{M_{p}}\right]^{\top} \in \mathbb{R}^{PM_{p}}$ with $P \geq 1$, $m_p = 1, \ldots, M_P$ and $G_{i} \in \mathbb{R}^{P\times P}$ is a known control input matrix. The control signal node is $u_{i} \in \mathbb{R}^{P}$ and the output vector is $y_{i} \in \mathbb{R}^{P}$ with $i = 1,\ldots,N$. It should be noted that the system has $p = 1, \ldots, P$ and $P$ is number of control inputs and it equals to number of regulated outputs. The nonlinear dynamics $f_{i}: \mathbb{R}^{P\times M_{p}} \rightarrow \mathbb{R}^{P}$ is unknown vector and Lipschitz. The global dynamics of equation \eqref{eq:equa1} can be described by
\begin{equation}
\label{eq:equa2}
\begin{aligned}
\dot{x}^{1} & = x^{2}\\
\dot{x}^{2} & = x^{3}\\
& \vdots\\
\dot{x}^{M_{p}} & = f\left(x\right) + Gu\\
y & = x^{1}
\end{aligned}
\end{equation}
where $x^{m_{p}} = \left[x_{1}^{m_{p}},\ldots,x_{N}^{m_{p}}\right]^{\top}\in \mathbb{R}^{PN}$, $u = \left[u_1^{\top},\ldots,u_N^{\top}\right]^{\top}\in \mathbb{R}^{PN}$, $G ={\rm  diag}\{G_{i}\}\in \mathbb{R}^{PN\times PN}$, $y = \left[y_1,\ldots,y_N\right]^{\top}\in \mathbb{R}^{P N}$,  $i = 1,\ldots,N$ and $f\left(x\right) = \left[f_1\left(x_1\right),\ldots,f_N\left(x_N\right)\right]^{\top}\in \mathbb{R}^{PN}$. The leader state vector can be time-varying and is noted $x_{0}$. It can be considered as an exosystem defining the desired consensus trajectory. Let's define the leader dynamics by
\begin{equation}
\label{eq:equa3}
\begin{aligned}
\dot{x}_{0}^{1} & = x_{0}^{2}\\
\dot{x}_{0}^{2} & = x_{0}^{3}\\
& \vdots\\
\dot{x}_{0}^{M_{p}} & = f_0\left(t,x_{0}\right)\\
y_{0} & = x_{0}^{1}
\end{aligned}
\end{equation}
where $x_{0}^{m_{p}} \in \mathbb{R}^P$ is the $m_p$-th state variable of the leader where $x_{0} = \left[x_{0}^{1}, \ldots, x_{0}^{M_P}\right]^{\top} \in \mathbb{R}^{PM_P}$ and the leader nonlinear function vector $f_0:\left[0,\infty\right) \times \mathbb{R}^{P\times M_{p}} \rightarrow \mathbb{R}^{P}$ is piecewise continuous in $t$ and locally Lipschitz. The disagreement variable for node $i$ is $\delta_{i}^{1} = x_{i}^{1} - x_{0}^{1}$ and the global disagreement order is  
\begin{equation}
\label{eq:equa4}
\begin{aligned}
& \gamma^{M_{p}} = x^{1} - \underline{x}_{0}^{1}
\end{aligned}
\end{equation}
where $\gamma^{M_{p}} = \left[\gamma_{1}^{1}, \ldots, \gamma_{N}^{1}\right]^{\top} \in \mathbb{R}^{PN} $, $\underline{x}_{0}^{1} = \left[x_{0}^{1}, \ldots, x_{0}^{1}\right]^{\top} \in \mathbb{R}^{PN}$. In this paper, local distributed state information is assumed to be known on the communication graph for $i$-th node and the only given information of the neighborhood synchronization as in \cite{li_pinning_2004,khoo_robust_2009} is  
\begin{equation}
\label{eq:equa5}
\begin{aligned}
& e_{i} = \sum_{j \in N_{i}} a_{i,j}\left(x_{i}^{1}-x_{j}^{1}\right)+b_{i,i}\left(x_{i}^{1}-x_{0}^{1}\right),
\end{aligned}
\end{equation}
where $e_{i} = \left[e_{i}^{1}, \ldots, e_{i}^{P}\right]^{\top} \in \mathbb{R}^{P}$, $a_{i,j}\geq 0$ and $a_{i,j} > 0$ in case of agent $i$ is directed to agent $j$, $b_{i}\geq 0$ and $b_{i}> 0$ for one or more agents $i$ are directed to the leader. $e^{p} = \left[e_{1}^{p}, \ldots, e_{N}^{p}\right]^{\top} \in \mathbb{R}^{PN}$, $p=1,\ldots,P$ and $B ={\rm  diag}\left\{b_{i}\right\}\in \mathbb{R}^{N \times N}$. The Global error dynamics for SISO system can be driven from \eqref{eq:equa5} to be
\begin{equation}
\label{eq:equa6}
\begin{aligned}
e &  = -\left(L+B\right)\left(\underline{x}_{0}^{1}-x^1\right)\\
& = \left(L+B\right)\left(x^1 - \underline{x}_{0}^{1}\right)
\end{aligned}
\end{equation}
Thereby, the error dynamics of \eqref{eq:equa6} can be written in the global form such as
\begin{equation}
\label{eq:equa7}
\begin{aligned}
\dot{e}^{1} & = e^{2}\\
\dot{e}^{2} & = e^{3}\\
& \vdots\\
\dot{e}^{M_{p}} & = \left(L+B\right)\left(f\left(x\right) + Gu - \underline{f}_0\right)\\
\end{aligned}
\end{equation}
Notice that $\underline{f}_0 = \left[f_0\left(t,x_{0}\right), \ldots, f_0\left(t,x_{0}\right)\right]^{\top} \in \mathbb{R}^{N}$. The proof of equation \eqref{eq:equa7} can be found in \cite{lewis_cooperative_2013}.\\

\begin{rem}
	Global error dynamics of \eqref{eq:equa6} for MIMO systems in case of $P>1$ becomes
	\begin{equation}
	\label{eq:equa8}
	\begin{aligned}
	e &  = -\left(\left(L+B\right)\otimes \mathbb{I}_P\right)\left(\underline{x}_{0}^{1}-x^1\right)\\
	& = \left(\left(L+B\right)\otimes \mathbb{I}_P\right)\left(x^1 - \underline{x}_{0}^{1}\right)
	\end{aligned}
	\end{equation}
	similarly, equation \eqref{eq:equa7} for the MIMO case can be written as
	\begin{equation}
	\label{eq:equa9}
	\begin{aligned}
	\dot{e}^{1} & = e^{2}\\
	\dot{e}^{2} & = e^{3}\\
	& \vdots\\
	\dot{e}^{M_{p}} & = \left(\left(L+B\right)\otimes \mathbb{I}_P\right)\left(f\left(x\right) + Gu - \underline{f}_0\right)\\
	\end{aligned}
	\end{equation}
	with $\otimes$ is the Kronecker product and $\mathbb{I}_P \in \mathbb{R}^{P\times P}$ is the identity matrix.
\end{rem} 

\begin{rem}
	The networked graph is strongly connected. Therefore, if there is one or more nodes $i$, $i=1,\ldots, N$ such that $b_{i} \neq 0$, then the matrix $\left(L+B\right)$ is an irreducible diagonally dominant M-matrix. Thus, it is nonsingular \cite{Qu2009}. 
\end{rem}

for the case of graph is strongly connected, $B \neq 0$ and the $\left\Vert e_0 \right\Vert$ is
\begin{equation}
\label{eq:equa10}
\begin{aligned}
& \left\Vert e_0 \right\Vert \leq \frac{\left\Vert e \right\Vert}{\underline{\sigma}\left(L+B\right)}
\end{aligned}
\end{equation}
such that $\underline{\sigma}\left(L+B\right)$ denotes the minimum singular value of matrix $L+B$.\\
%\begin{rem}
%In the case of multi-agent systems, Equation (\ref{eq:equa6}) reflects the coupling that has been creating through the synchronization between the different states of each agent. Thus, these interactions no longer allows to guarantee that the error dynamics of each agent will be confined within the desired performance functions just based on knowledge of the sign $e_{i}\left(0\right)$, $i=1,...N$.    
%\end{rem}

\subsection{Prescribed Performance}
\begin{figure*}[ht]
	\centering
	\includegraphics[scale=0.4]{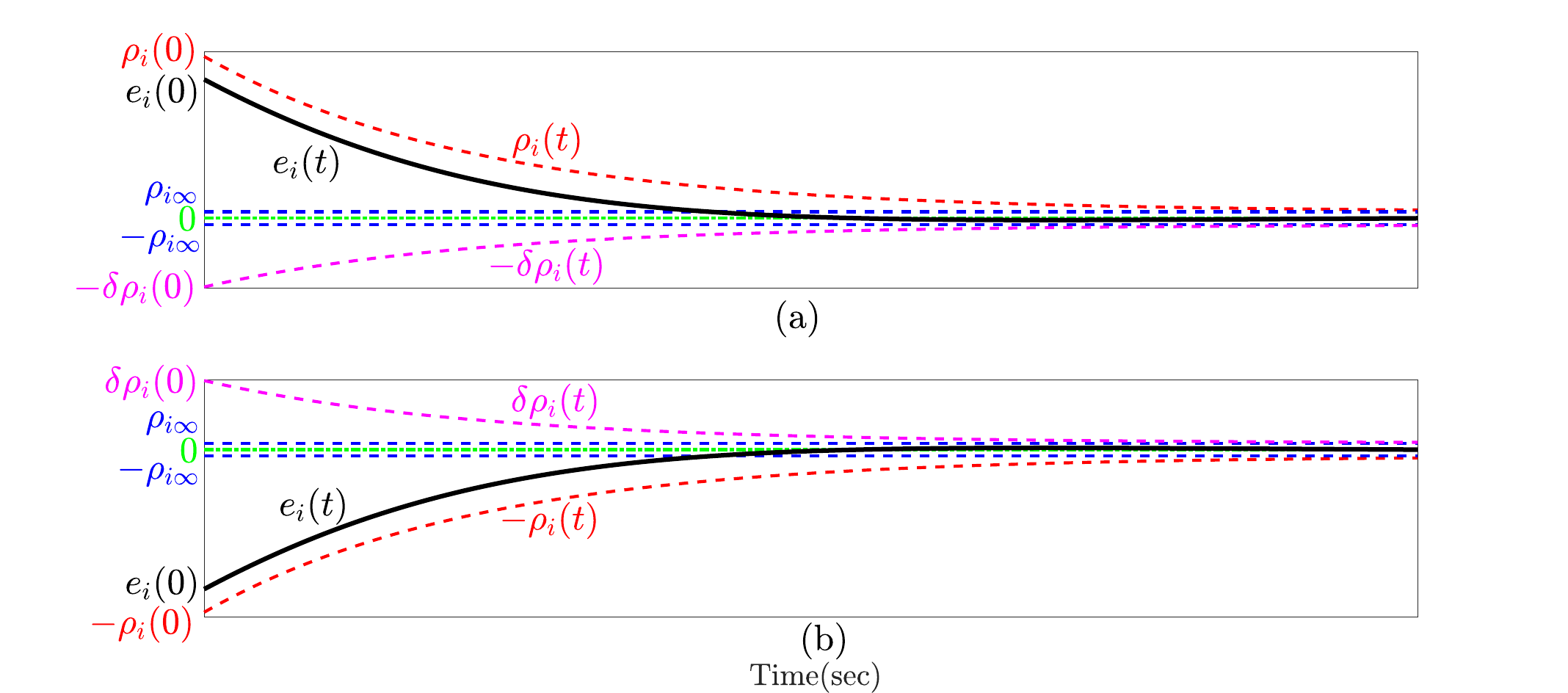}
	\caption{ Schematic representation of error trajectory with prescribed performance
		(a) Schematic illustration of prescribed performance of \eqref{eq:equa12}; (b) Graphical illustration of prescribed performance of \eqref{eq:equa13}.}
	\label{fig:fig1}
\end{figure*}

The objective of this subsection is to introduce the prescribed performance function (PPF) into the control algorithm. PPF is a time function enables the tracking error $e\left(t\right)$ to start within a known large set and reduce in a systematic manner to a known narrow set \cite{bechlioulis_robust_2008} and \cite{Hashim2017adaptive}. Providing smooth tracking response with allocated properties and improving the control signal range are classified as distinguished features of the control algorithm with PPF.\\
Consider the performance function of a single agent system with $\rho\left(t\right)$ is a smooth function includes the error component $e\left(t\right)$ such as $\rho\left(t\right): \mathbb{R}_{+} \to \mathbb{R}_{+}$ is a decreasing positive function $\lim\limits_{t \to \infty}\rho\left(t\right)=\rho_{\infty}>0$ where $\rho_{\infty}>0$ is a constant and refers to the smaller set upper bound. Now, the general PPF of \eqref{eq:equa5} can be described as
\begin{equation}
\label{eq:equa11}
\rho_{i}^{p}\left(t\right)=\left(\rho_{i,0}^{p} - \rho_{i,\infty}^{p}\right)\exp\left(-\ell_{i}^{p}t\right)+\rho_{i,\infty}^{p}
\end{equation}

where $\rho_{i}\left(t\right) = \left[\rho_{i}^{1}, \ldots, \rho_{i}^{P}\right]^{\top} \in \mathbb{R}^{P}$, also, $\rho_{i,0}^{p}$, $\rho_{i,\infty}^{p}$ and $\ell_{i}^{p}$ are appropriately defined positive constants with $p=1,\ldots,P$. $\rho_{i,0}^{p}$ and $\rho_{i,\infty}^{p}$ are positive constants define initial and final upper bounds of the predefined sets. The prescribed properties of the control function should guarantee the following properties:
%      \begin{equation}
%         \label{eq:equa12}
%         -\delta_{i}^{p}\rho_{i}^{p}\left(t\right)<e_{i}^{p}\left(t\right)<\rho_{i}^{p}\left(t\right),\hspace{1pt} {\rm{if}} \: e_{i}^{p}\left(0\right)>0
%       \end{equation}
%       \begin{equation}
%         \label{eq:equa13}
%         -\rho_{i}^{p}\left(t\right)<e_{i}^{p}\left(t\right)<\delta_{i}^{p}\rho_{i}^{p}\left(t\right),\hspace{1pt} {\rm{if}} \: e_{i}^{p}\left(0\right)<0
%       \end{equation}
\begin{align}
-\delta_{i}^{p}\rho_{i}^{p}\left(t\right)<e_{i}^{p}\left(t\right)<\rho_{i}^{p}\left(t\right)&,\hspace{2pt} {\rm{if}} \: e_{i}^{p}\left(0\right)>0\label{eq:equa12}\\
-\rho_{i}^{p}\left(t\right)<e_{i}^{p}\left(t\right)<\delta_{i}^{p}\rho_{i}^{p}\left(t\right)&,\hspace{2pt} {\rm{if}} \: e_{i}^{p}\left(0\right)<0\label{eq:equa13}
\end{align}       

for all $t \geq 0 $ and $ 0 \leq \delta_{i}^{p} \leq 1 $, $i=1, \ldots,N$ and $p=1, \ldots,P$. The control algorithm should consider the interactions between agents' dynamics which may lead to instability. The systematic convergence of the tracking error $e_{i}^{p}\left(t\right)$ between the constraint bounds $\rho_{i}^{p}\left(t\right)$ and $-\delta_{i}^{p}\rho_{i}^{p}\left(t\right)$ or $-\rho_{i}^{p}\left(t\right)$ and $\delta_{i}^{p}\rho_{i}^{p}\left(t\right)$ should obey the transient trajectory of these foregoing bounds as revealed in Figure \ref{fig:fig1}. In fact, Figure \ref{fig:fig1} illustrates the full idea of prescribed performance such that the error will be tracked systematically from a pre-defined bigger set to a given smaller set.
A transformed error will be defined to drive the error dynamics from constrained bounds in \eqref{eq:equa12} and \eqref{eq:equa13} into an unconstrained one as follows 
\begin{equation}
\label{eq:equa14}
\varepsilon_{i}^{p}=\Upsilon\left(\frac{e_{i}^{p}\left(t\right)}{\rho_{i}^{p}\left(t\right)}\right)
\end{equation}
or equivalently,
\begin{equation}
\label{eq:equa15}
e_{i}^{p}\left(t\right)=\rho_{i}^{p}\left(t\right)\mathcal{F}\left(\varepsilon_{i}^{p}\right)
\end{equation}
where $\varepsilon_{i}^{p}$, $\mathcal{F}\left(\cdot\right)$ and $\Upsilon^{-1}\left(\cdot\right)$ are smooth functions, $i=1,2,\ldots,N$. For simplification, let us denote $x:=x\left(t\right)$, $\rho:=\rho\left(t\right)$, $e:=e\left(t\right)$ and $\varepsilon:=\varepsilon\left(t\right)$.  $\mathcal{F}\left(\cdot\right)=\Upsilon^{-1}\left(\cdot\right)$ and $\mathcal{F}\left(\cdot\right)$ satisfy the following properties:
\begin{enumerate}
	\item $\mathcal{F}\left(\varepsilon_{i}^{p}\right)$ is smooth and strictly increasing.
	\item $-\underline{\delta}_{i}^{p}<\mathcal{F}\left(\varepsilon_{i}^{p}\right)<\bar{\delta}_{i}^{p},\hspace{5pt} {\rm if} \hspace{2pt} e_{i}^{p}\left(0\right) \geq 0$\\
	$-\bar{\delta}_{i}^{p}<\mathcal{F}\left(\varepsilon_{i}^{p}\right)<\underline{\delta}_{i}^{p},\hspace{5pt} {\rm if} \hspace{2pt} e_{i}^{p}\left(0\right)<0$
	\item
	$\left.
	\begin{aligned}
	{\rm lim}_{\varepsilon_{i}^{p} \rightarrow -\infty}\mathcal{F}\left(\varepsilon_{i}^{p}\right)=-\underline{\delta}_{i}^{p}\\
	{\rm lim}_{\varepsilon_{i}^{p} \rightarrow +\infty}\mathcal{F}\left(\varepsilon_{i}^{p}\right)=\bar{\delta}_{i}^{p}
	\end{aligned}
	\right\}
	\quad {\rm if} \hspace{5pt} e_{i}^{p}\left(0\right)\geq 0$\\
	$\left.
	\begin{aligned}
	{\rm lim}_{\varepsilon_{i}^{p} \rightarrow -\infty}\mathcal{F}\left(\varepsilon_{i}^{p}\right)=-\bar{\delta}_{i}^{p}\\
	{\rm lim}_{\varepsilon_{i}^{p} \rightarrow +\infty}\mathcal{F}\left(\varepsilon_{i}^{p}\right)=\underline{\delta}_{i}^{p}
	\end{aligned}
	\right\}
	\quad {\rm if} \hspace{5pt} e_{i}^{p}\left(0\right)< 0$
\end{enumerate}

for $\underline{\delta}_{i}^{p}, \bar{\delta}_{i}^{p}\in \mathbb{R}_{+}$ are known constants. These constants should be defined to satisfy
\small{
	\begin{equation}
	\label{eq:equa16}
	\mathcal{F}\left(\varepsilon_{i}^{p}\right)=
	\left\{
	\begin{aligned}
	\frac{\bar{\delta}_{i}^{p}\exp\left(\varepsilon_{i}^{p}\right)-\underline{\delta}_{i}^{p}\exp\left(-\varepsilon_{i}^{p}\right)}{\exp\left(\varepsilon_{i}^{p}\right)+\exp\left(-\varepsilon_{i}^{p}\right)},& \quad \bar{\delta}_{i}^{p}> \underline{\delta}_{i}^{p} \hspace{3pt} {\rm if} \hspace{3pt} e_{i}^{p}\left(0\right)\geq 0\\
	\frac{\bar{\delta}_{i}^{p}\exp\left(\varepsilon_{i}^{p}\right)-\underline{\delta}_{i}^{p}\exp\left(-\varepsilon_{i}^{p}\right)}{\exp\left(\varepsilon_{i}^{p}\right)+\exp\left(-\varepsilon_{i}^{p}\right)},& \quad \underline{\delta}_{i}^{p} > \bar{\delta}_{i}^{p} \hspace{3pt} {\rm if} \hspace{3pt} e_{i}^{p}\left(0\right) < 0\\
	\end{aligned}
	\right.
	\end{equation}
}
Now, consider the smooth function
\begin{equation}
\label{eq:equa17}
\mathcal{F}\left(\varepsilon_{i}^{p}\right)=
\begin{aligned}
& \frac{\bar{\delta}_{i}^{p}\exp\left(\varepsilon_{i}^{p}\right)-\underline{\delta}_{i}^{p}\exp\left(-\varepsilon_{i}^{p}\right)}{\exp\left(\varepsilon_{i}^{p}\right)+\exp\left(-\varepsilon_{i}^{p}\right)}
\end{aligned}
\end{equation}
and the transformed error
\begin{equation}
\label{eq:equa18}
\begin{aligned}
\varepsilon_{i}^{p} = & \mathcal{F}^{-1}\left(e_{i}^{p}/\rho_{i}^{p}\right)\\
= & \frac{1}{2}\left\{ \begin{aligned} & {\rm ln} \frac{\underline{\delta}_{i}^{p}+e_{i}^{p}/\rho_{i}^{p}}{\bar{\delta}_{i}^{p}-e_{i}^{p}/\rho_{i}^{p}}, & \bar{\delta}_{i}^{p} > \underline{\delta}_{i}^{p} \hspace{3pt} {\rm if} \hspace{3pt} e_{i}^{p}\left(0\right)\geq 0 \\
& {\rm ln} \frac{\underline{\delta}_{i}^{p}+e_{i}/\rho_{i}^{p}}{\bar{\delta}_{i}^{p}-e_{i}^{p}/\rho_{i}^{p}} ,& \underline{\delta}_{i}^{p} > \bar{\delta}_{i}^{p} \hspace{3pt} {\rm if} \hspace{3pt} e_{i}^{p}\left(0\right)< 0\\ \end{aligned}  \right.
\end{aligned}
\end{equation}

%     \begin{rem} The primary role of $\xi$ is to make $erf(\xi e)$ as close as possible to $sign(e)$. Ideally $\xi$ is selected as big as possible. For instance, $|erf(\xi e)|\approxeq 1$ when $ |e|>\Delta=\frac{2}{\xi}$. Therefore, if $\xi=200$ then $|erf(e)|\approxeq 1$ when $|e|>0.01$.  However, while the derivative is smooth the more one selects a big $\xi$ the more there is a risk of chattering. 
%     \end{rem}  
Therefore, the derivative of the transformed error in \eqref{eq:equa18} will be
\begin{equation}
\label{eq:equa19}
\dot{\varepsilon}_{i}^{p} = \frac{1}{2\rho_{i}^{p}}\left(\frac{1}{\underline{\delta}_{i}^{p}+e_{i}^{p}/\rho_{i}^{p}} + \frac{1}{\bar{\delta}_{i}^{p} - e_{i}^{p}/\rho_{i}^{p}}\right)\left(\dot{e}_{i}^{p} - \frac{e_{i}^{p}\dot{\rho}_{i}^{p}}{\rho_{i}^{p}}\right)
\end{equation}
where $\varepsilon_{i} \in \mathbb{R}^{P}$ and from \eqref{eq:equa19}, we define new variable $r_{i}^{p}$ such as
\begin{equation}
\begin{split}
\label{eq:equa20}
r_{i}^{p} & = \frac{1}{2\rho_{i}^{p}}\frac{\partial \mathcal{F}^{-1}\left(e_{i}^{p}/\rho_{i}^{p}\right)}{\partial\left(e_{i}^{p}/\rho_{i}^{p}\right)}\\
& = \frac{1}{2\rho_{i}^{p}}\left(\frac{1}{\underline{\delta}_{i}^{p}+e_{i}^{p}/\rho_{i}^{p}} + \frac{1}{\bar{\delta}_{i}^{p} - e_{i}^{p}/\rho_{i}^{p}}\right)
\end{split}
\end{equation}
For further explanations, we define a new component ${\bf E}_{i} \in \mathbb{R}^{P}$ such that ${\bf E}_{i}$ is a metric error that can be described as
\begin{equation}
\label{eq:equa21}
{\bf E}_{i} = \left(\frac{d}{dt} + \lambda_{i}^{m_{p}}\right)^{M_{p}-1}\varepsilon_{i}^{1}, \hspace{5pt} i = 1,\ldots,N, \hspace{5pt} m_p = 1,\ldots,M_{p}
\end{equation}
where $\lambda_{i}^{m_{p}}$ is a positive constant, alternatively, \eqref{eq:equa21} is equivalent to
\begin{equation}
\label{eq:equa22}
{\bf E}_{i} = \varepsilon_{i}^{M_{p}} + \lambda_{i}^{M_{p}-1}\varepsilon_{i}^{M_{p}-1} + \cdots + \lambda_{i}^{1}\varepsilon_{i}^{1}
\end{equation}
one can write the global form of \eqref{eq:equa22} as
\begin{equation}
{\bf E} = \varepsilon^{M_{p}} + \lambda^{M_{p}-1}\varepsilon^{M_{p}-1} + \cdots + \lambda^{1}\varepsilon^{1}
\end{equation}
where $\varepsilon^{m_p} = \left[\varepsilon_{1}^{m_p}, \ldots, \varepsilon_{N}^{m_p}\right]^{\top}$, $m_p = 1, \ldots, M_P$.
Let's define
\begin{align}
\Phi_{1} & =\left[\varepsilon^{1},\varepsilon^{2},\ldots,\varepsilon^{Mp-1}\right]^{\top}\label{eq:Case1}\\
\Phi_{2} & =\dot{\Phi}_{1}=\left[\varepsilon^{2},\varepsilon^{3},\ldots,\varepsilon^{Mp}\right]^{\top}\label{eq:Case2}\\
l & =\left[0,0,\ldots,0,1\right]^{\top}\in\mathbb{R}^{M_{p}-1}\nonumber 
\end{align}
and
\begin{equation*}
\begin{aligned}
& \Lambda = \begin{bmatrix}
0 & 1 & 0 & \cdots & 0 & 0\\
0 & 0 & 1 & \cdots & 0 & 0\\
\vdots & \vdots & \vdots & \ddots & \vdots & \vdots\\
0 & 0 & 0 & \cdots & 0 & 1\\
-\lambda^1 & -\lambda^2 & -\lambda^3 & \cdots & \lambda^{M_{p}-2} & -\lambda^{M_{p}-1}\\
\end{bmatrix}
\end{aligned}
\end{equation*}
such that $\Lambda$ is Hurwitz , hence, from \eqref{eq:Case1} and \eqref{eq:Case2} one can say
\begin{equation}
\begin{aligned}
\label{eq:equa24}
& \Phi_{2} = \Phi_{1}\Lambda^{\top} + {\bf E}l^{\top}\\
\end{aligned}
\end{equation}
and
\begin{equation}
\begin{aligned}
\label{eq:equa25}
& \Lambda^{\top}M + M\Lambda = -\beta \mathbb{I}_{M_{p}-1}
\end{aligned}
\end{equation}
where $\beta$ is a positive constant, $M>0$ and $\mathbb{I}_{M_{p}-1} \in \mathbb{R}^{\left(M_{p}-1\right)\times \left(M_{p}-1\right)}$ is the identity matrix. Consider each of  \eqref{eq:equa7} and \eqref{eq:equa18}, the derivative of the metric error in \eqref{eq:equa21} with respect to time is given by
\begin{equation}
\begin{aligned}
\label{eq:equa26}
\dot{{\bf E} }_{i} = & \sum\limits_{j=1}^{M_{p}-1} \begin{bmatrix} M_{p}-1\\ j \end{bmatrix} \lambda_{i}^{j}\varepsilon_{i}^{M_{p}-j} + \varepsilon_{i}^{M_{p}}\\
\end{aligned} 
\end{equation}
with $\bar{\lambda} = \left[\lambda^1, \ldots, \lambda^{M_{p}-1}\right]^{\top} \in \mathbb{R}^{PN}$, $p = 1, \ldots, P$, and equation \eqref{eq:equa26} can be written in the global form for SISO systems as
\begin{equation}
\begin{aligned}
\label{eq:equa27}
\dot{{\bf E} } & =  \varepsilon^{M_{p}+1} + \Phi_{2}\bar{\lambda}\\
& = R\left(L+B\right)\left(f\left(x\right) + Gu - \underline{f}_0\right) + \Delta + \Phi_{2}\bar{\lambda}\\
\end{aligned} 
\end{equation}
and for MIMO case
\begin{equation}
\begin{aligned}
\label{eq:equa28}
\dot{{\bf E} } & = R\left(\left(L+B\right)\otimes \mathbb{I}_{P}\right)\left(f\left(x\right) + Gu - \underline{f}_0\right) + \Delta + \Phi_{2}\bar{\lambda}\\
\end{aligned} 
\end{equation}
where ${\bf E}=\left[{\bf E}_{1},\ldots,{\bf E}_{N}\right]^{\top}\in \mathbb{R}^{PN}$, $\varepsilon^{M_{p}+1} = e^{M_{p}+1} + \Delta$ and $\Delta$ is the function of higher orders of $\rho_{i}^{p}$, $r_{i}^{p}$. It should be remarked that the higher orders of $\rho_{i}^{p}$, $r_{i}^{p}$ are vanishing components with time which by the way lead to $\Delta = 0$ as $t \rightarrow \infty$. Also, $R ={\rm  diag}\left\{\Omega_{i}\right\} \in \mathbb{R}^{PN \times PN}$ and
%      $
%       \dot{\rho}_{i}}{\rho_1},\ldots,\frac{\dot{\rho}_N}{\rho_N}]\in \mathbb{R}^N $ and $R ={\rm  diag}\{r_{i}\}, i = 1,\ldots,N$ or  
\begin{equation*}
\Omega_{i} =
\begin{pmatrix}
\frac{1}{2\rho_{i}^{1}}\frac{\partial \mathcal{F}^{-1}\left(e_{i}^{1}/\rho_{i}^{1}\right)}{\partial\left(e_{i}^{1}/\rho_{i}^{1}\right)} & \cdots & 0\\
\vdots & \ddots & \vdots\\
0 & \cdots & \frac{1}{2\rho_{i}^{P}} \frac{\partial \mathcal{F}^{-1}\left(e_{i}^{P}/\rho_{i}^{P}\right)}{\partial\left(e_{i}^{P}/\rho_{i}^{P}\right)}\\
\end{pmatrix}
\end{equation*}
Note that $\Omega_{i}$ is a decreasing diagonal matrix with $\Omega_{i}>0$. Let's recall the following definitions (see \cite{das_distributed_2010})
\begin{defn} \label{Def1}
	The global neighborhood error $e\left(t\right)\in \mathbb{R}^{PN} $ is uniformly ultimately bounded (UUB) if there exists a compact set $\Psi \subset \mathbb{R}^{PN} $ so that $\forall e\left(t_{0}\right) \in \Psi$ there exists a bound $B$ and a time $t_f(B,e\left(t_{0}\right))$, both independent of $t_{0} \geq 0$, such that $\left\Vert e\left(t\right) \right\Vert \leq B$ so that $\forall t > t_{0}+t_f$.
\end{defn}
\begin{defn} \label{Def2}
	The control node trajectory $x_{0}\left(t\right)$ given by \eqref{eq:equa1} is cooperative UUB with respect to solutions of node dynamics \eqref{eq:equa3} if there exists a compact set $\Psi \subset \mathbb{R}^{PN} $ so that $\forall \left(x_{i}\left(t_{0}\right)-x_{0}\left(t_{0}\right)\right) \in \Psi$, there exist a bound $B$ and a time $t_m\left(B,\left(x\left(t_{0}\right) - x_{0}\left(t_{0}\right)\right)\right)$, both independent of $t_{0}\geq 0$, such that $\left\Vert \gamma^{m_{p}} \right\Vert \leq B$, $\forall i$, $m_{p}=1,\ldots,M_{p}$ and $\forall t > t_{0}+t_m$.
\end{defn}
\section{Neural Approximation and Distributed Control in Prescribed Performance} \label{Sec4}
\subsection{Neural Approximations}
Neural network with linear weights are used to approximate the unknown nonlinear dynamics of local agents in \eqref{eq:equa1} as   
\begin{equation}
\label{eq:equa29}
f_{i}\left(x_{i}\right) = W_{i}^{\top}\phi_{i}\left(x_{i}\right) + \alpha_{i}
\end{equation}
where $\phi_{i}\left(x_{i}\right) \in \mathbb{R}^{v_{i} \times 1}$ and $v_{i}$ is a sufficient number of neurons at each node, $W_{i} \in \mathbb{R}^{v_{i} \times P}$ and $\alpha_{i} \in \mathbb{R}^{P \times 1}$ is the approximated error vector. It should be remarked that based on \cite{hornik_multilayer_1989,lewis_neural_1998}, nonlinearities could be approximated via variety of sets such as radial basis functions \cite{poggio_regularization_1990}, sigmoid functions \cite{cotter_stone-weierstrass_1989}, etc.\\
The objective of this work is to use the available information to track local performance behavior of each node and to compensate unknown nonlinearities. Thereby, the unknown nonlinearities of the local nodes can be approximated by
\begin{equation}
\label{eq:equa30}
\hat{f}_{i}\left(x_{i}\right) = \hat{W}_{i}^{\top}\phi_{i}\left(x_{i}\right)
\end{equation}
where $\hat{W}_{i} \in \mathbb{R}^{v_{i} \times P}$ and $\hat{f}_{i}\left(x_{i}\right)\in \mathbb{R}^{P}$ approximate the component $f_{i}\left(x_{i}\right)$. The description of the global synchronization of the graph $\mathcal{G}$ could be defined by
\begin{equation}
\label{eq:equa31}
f\left(x\right) = W^{\top}\phi\left(x\right) + \alpha
\end{equation}
where $\phi\left(x\right)=\left[\phi_1\left(x_1\right),\ldots,\phi_N\left(x_N\right)\right]^{\top}$ $,i=1,\ldots,N$, $W={\rm diag}\left\{W_{i}\right\}$, $\alpha=\left[\alpha_1,\ldots,\alpha_N\right]^{\top}$ and the global estimate of $f\left(x\right)$ is
\begin{equation}
\label{eq:equa32}
\hat{f}\left(x\right) = \hat{W}^{\top}\phi\left(x\right)
\end{equation}
with $\hat{W}^{\top}={\rm diag}\left\{\hat{W}_{i}\right\}$. The error between true and estimated nonlinearities is defined as
\begin{equation}
\label{eq:equa33}
\tilde{f}\left(x\right) = f\left(x\right) - \hat{f}\left(x\right) =\tilde{W}^{\top}\phi\left(x\right) + \alpha
\end{equation}
such as $\tilde{W}= W - \hat{W}$.

\subsection{Neuro-Adaptive Control Design with PPF of Distributed Agents}
There are several assumptions should be considered \cite{zhang_adaptive_2012}
\begin{assum}\label{assum1}
	\begin{enumerate}
		\item Neural network (NN) weights are bounded but otherwise, they are unknown such as $\left\Vert W \right\Vert \leq W_M$  with $W_M$ is a fixed bound. 
		\item Leader states $\left\Vert x_{0} \right\Vert \leq X_0$ are bounded with $X_0$ is a finite bound.
		\item The unknown nonlinear dynamics associated to the leader is bounded by $\left\Vert \underline{f}\left(x_{0},t\right) \right\Vert \leq F_M$.
		\item The activation function is finite such as $\left\Vert \phi \right\Vert \leq \phi_M$.
	\end{enumerate}
\end{assum}
\begin{lem}
	\label{lemma1} {\bf \cite{Qu2009}}\\
	Consider $L$ is an irreducible matrix with $\left(L+B\right)$ is a nonsingular matrix such as and $B \neq 0 $, hence there is
	\begin{equation}
	\label{eq:equa34}
	q = [q_{1},\ldots,q_{N}]^{\top}=\left(L+B\right)^{-1}\cdot\underline{1}
	\end{equation}
	\begin{equation}
	\label{eq:equa35}
	\mathcal{M} ={\rm  diag}\left\{m_{i}\right\}={\rm diag}\left\{1/q_{i}\right\}
	\end{equation}
	Then, $\mathcal{M} > 0$ and the matrix $\mathcal{Q}$ defined as
	\begin{equation}
	\label{eq:equa36}
	\mathcal{Q} =\mathcal{M} \left(L+B\right)+\left(L+B\right)^{\top}\mathcal{M}
	\end{equation}
\end{lem}
\begin{rem}
	For simplification, the control signal and the estimated weights of NN throughout the research paper will be developed for SISO systems. In case of MIMO systems, these developments can be easily modified by introducing the Kronecker product as will be illustrated later. Two different simulation for SISO and MIMO systems will be presented.
\end{rem}
Consider a control signal of each local nodes be
\begin{equation}
\begin{aligned}
\label{eq:equa37}
u_{i} = & B_{i}^{-1}\bigg(-c{\bf E}_{i} - \hat{W_{i}}^{\top}\phi_{i}\left(x_{i}\right)\\
& \hspace{19pt} - \left(d_{i}+b_{i}\right)^{-1}\Omega_{i}^{-1}( \lambda_{i}^{M_{p}-1}\varepsilon_{i}^{M_{p}} + \cdots + \lambda_{i}^{1}\varepsilon_{i}^{2})\bigg)
\end{aligned}
\end{equation}
where $c$ be a positive control gain and ${\bf \underline{1}}_{N\times 1} = \left[1,\ldots,1\right]^{\top}\in \mathbb{R}^{N \times 1}$. The control input is defined by
\begin{equation}
\label{eq:equa38}
u = G^{-1}\left(-c{\bf E} - \hat{W}^{\top}\phi\left(x\right) - \left(D+B\right)^{-1}R^{-1}\Phi_{2}\bar{\lambda}\right)
\end{equation}
Then, the estimated weights of NN is defined by
\begin{equation}
\label{eq:equa39}
\dot{\hat{W}}_{i} = F_{i}\phi_{i}{\bf E}_{i}^{\top}m_{i}\Omega_{i}\left(d_{i}+b_{i}\right) - kF_{i}\hat{W}_{i}
\end{equation}
with $F_{i} \in \mathbb{R}^{v_{i} \times v_{i}}$, $F_{i} = \Pi_{i}\mathbb{I}_{v_{i}}$ and $\Pi_{i} >0$ are positive gains and $k > 0$ is a scalar gain. $c$ and $k$ should be selected to satisfy \eqref{eq:equa40}.
\begin{thm}\label{theom1}
	Consider the distributed system in \eqref{eq:equa1} and the leader dynamics in \eqref{eq:equa3}. If Assumption \ref{assum1} holds and the distributed control is as in \eqref{eq:equa38} and the NN tuning law as in \eqref{eq:equa39}, then the control variable $c$ should satisfy
	\begin{equation}
	\label{eq:equa40}
	c> \frac{1}{\underline{\sigma}\left(\mathcal{Q}\right)\underline{\sigma}\left(R\right)}\left(\frac{\gamma^2}{k}+\frac{2}{\beta}g^2 + \nu\right)
	\end{equation}
	$\gamma = -\frac{1}{2}\Phi\bar{\sigma}\left(\mathcal{M}\right)\bar{\sigma}\left(R\right)\bar{\sigma}\left(A\right)$, $g = -\frac{1}{2}\left(\bar{\sigma}\left(M\right)+\frac{\bar{\sigma}\left(\mathcal{M}\right)\bar{\sigma}\left(A\right)}{\underline{\sigma}\left(D+B\right)}\left\Vert \Lambda \right\Vert_F\left\Vert \bar{\lambda}\right\Vert\right)$ and $\nu = \frac{\bar{\sigma}\left(\mathcal{M}\right)\bar{\sigma}\left(A\right)}{\underline{\sigma}\left(D+B\right)}\left\Vert \bar{\lambda}\right\Vert$, where $M$ was defined in \eqref{eq:equa25} for $\beta >0$. Hence, the trajectory of $x_{0}\left(t\right)$ is uniformly ultimate bounded. Also, all nodes steer close to $x_{0}\left(t\right)$ for all $t\geq 0$.
\end{thm}

\textbf{Proof:}\\
Based on \eqref{eq:equa31}, equation \eqref{eq:equa7} becomes 
\begin{equation}
\dot{e}^{M_{p}}=\left(L+B\right)\left(W^{\top}\phi\left(x\right)+\alpha+Gu-\underline{f}\left(x_{0},t\right)\right) \label{eq:equa41}
\end{equation}
consider the result in \eqref{eq:equa32} and \eqref{eq:equa33}, one can write \eqref{eq:equa41} as 
\begin{equation}
\dot{e}=\left(L+B\right)\left(\tilde{W}^{\top}\phi\left(x\right)+\alpha-c{\bf E}-\left(D+B\right)^{-1}R^{-1}\Phi_{2}\bar{\lambda}-\underline{f}\left(x_{0},t\right)\right)\label{eq:equa42}
\end{equation}
and from \eqref{eq:equa27} and \eqref{eq:equa42}, the transformed error could be obtained as
\begin{equation}
\begin{aligned}
\dot{{\bf E} }= & R\left(L+B\right)\left(\tilde{W}^{\top}\phi\left(x\right)+\alpha-c{\bf E}-\left(D+B\right)^{-1}R^{-1}\Phi_{2}\bar{\lambda}-\underline{f}\left(x_{0},t\right)\right)\\
& +\Delta+\Phi_{2}\bar{\lambda}\label{eq:equa43}
\end{aligned}
\end{equation}
%\subsection{Lyapunov Functions of Distributed Control }
Consider the following Lyapunov candidate function 
\begin{equation}
\begin{aligned}\label{eq:equa44}V & =\frac{1}{2}{\bf E}^{\top}\mathcal{M}{\bf E}+\frac{1}{2}{\rm Tr}\left\{ \tilde{W}^{\top}F^{-1}\tilde{W}\right\} +\frac{1}{2}{\rm Tr}\left\{ \Phi_{1}M\Phi_{1}^{\top}\right\} \\
& =V_{1}+V_{2}+V_{3}
\end{aligned}
\end{equation}
with $V_{1}=\frac{1}{2}{\bf E}^{\top}\mathcal{M}{\bf E}$, $V_{2}=\frac{1}{2}{\rm Tr}\left\{ \tilde{W}^{\top}F^{-1}\tilde{W}\right\} $,
$V_{3}=\frac{1}{2}{\rm Tr}\left\{ \Phi_{1}M\Phi_{1}^{\top}\right\} $,
and $M>0$ is as defined in Lemma \ref{lemma1} and $F^{-1}={\rm diag}\{F_{i}^{-1}\}$ is a zero matrix with positive components in diagonal as in \eqref{eq:equa39}.
Then, $\dot{V_{1}}$ and $\dot{V_{2}}$ after substitution
of \eqref{eq:equa38} is 
\begin{equation}
\dot{V_{1}}+\dot{V_{2}}={\bf E}^{\top}\mathcal{M}\dot{{\bf E} }+{\rm Tr}\left\{ \tilde{W}^{\top}F^{-1}\dot{\tilde{W}}\right\} \label{eq:equa45}
\end{equation}
\begin{equation}
\begin{split}\dot{V_{1}}+\dot{V_{2}}= & {\bf E}^{\top}\mathcal{M}R\left(L+B\right)\Big(\tilde{W}^{\top}\phi\left(x\right)+\alpha-c{\bf E}\\
& -\left(D+B\right)^{-1}R^{-1}\Phi_{2}\bar{\lambda}-\underline{f}\left(x_{0},t\right)\Big) +{\bf E}^{\top}\mathcal{M}\left(\Delta+\Phi_{2}\bar{\lambda}\right)\\
& +{\rm Tr}\left\{ \tilde{W}^{\top}F^{-1}\dot{\tilde{W}}\right\} 
\end{split}
\label{eq:equa46}
\end{equation}
\begin{equation}
\begin{split}\dot{V_{1}}+\dot{V_{2}}= & -c{\bf E}^{\top} \mathcal{M}R\left(L+B\right){\bf E}+{\bf E}^{\top} \mathcal{M}R\left(D+B\right)\tilde{W}^{\top}\phi\left(x\right)\\
& -{\bf E}^{\top}\mathcal{M}RA\left(\tilde{W}^{\top}\phi\left(x\right)+\left(D+B\right)^{-1}R^{-1}\Phi_{2}\bar{\lambda}\right)\\
& +{\bf E}^{\top}\mathcal{M}R\left(L+B\right)\left(\alpha-\underline{f}\left(x_{0},t\right)\right)+{\bf E}^{\top} \mathcal{M}\Delta\\
& +{\rm Tr}\left\{ \tilde{W}^{\top}F^{-1}\dot{\tilde{W}}\right\} 
\end{split}
\label{eq:equa47}
\end{equation}
Note that $x^{\top}y={\rm Tr}\{yx^{\top}\}$, $\forall x,y\in\mathbb{R}^{N}$,
one can write $\dot{V_{1}}+\dot{V_{2}}$ as 
\begin{equation}
\begin{split}\dot{V_{1}}+\dot{V_{2}}= & -c{\bf E}^{\top}\mathcal{M}R\left(L+B\right){\bf E}+{\rm Tr}\left\{ \tilde{W}^{\top}\phi\left(x\right){\bf E}^{\top} \mathcal{M}R\left(D+B\right)\right\} \\
& -{\rm Tr}\left\{ \tilde{W}^{\top}\phi\left(x\right){\bf E}^{\top}\mathcal{M}RA\right\} +{\bf E}^{\top}\mathcal{M}RA\left(D+B\right)^{-1}R^{-1}\Phi_{2}\bar{\lambda}\\
& +{\bf E}^{\top}\mathcal{M}R\left(L+B\right)\left(\alpha-\underline{f}\left(x_{0},t\right)\right)+{\bf E}^{\top}\mathcal{M}\Delta\\
& +{\rm Tr}\left\{ \tilde{W}^{\top}F^{-1}\dot{\tilde{W}}\right\} 
\end{split}
\label{eq:equa48}
\end{equation}
substitute \eqref{eq:equa39} in \eqref{eq:equa48} considering $\dot{\tilde{W}}=\dot{W}-\dot{\hat{W}}=-\dot{\hat{W}}$ yields
\begin{equation}
\begin{split}\dot{V_{1}}+\dot{V_{2}}= & -c{\bf E}^{\top}\mathcal{M}R\left(L+B\right){\bf E}+{\rm Tr}\left\{ \tilde{W}^{\top}\phi\left(x\right){\bf E}^{\top}\mathcal{M}R\left(D+B\right)\right\} \\
& -{\rm Tr}\left\{ \tilde{W}^{\top}\phi\left(x\right){\bf E}^{\top} \mathcal{M}RA\right\} +{\bf E}^{\top}\mathcal{M}RA\left(D+B\right)^{-1}R^{-1}\Phi_{2}\bar{\lambda}\\
& +{\bf E}^{\top}\mathcal{M}R\left(L+B\right)\left(\alpha-\underline{f}\left(x_{0},t\right)\right)+{\bf E}^{\top}\mathcal{M}\Delta\\
& -{\rm Tr}\left\{ \tilde{W}^{\top}\phi\left(x\right){\bf E}^{\top}\mathcal{M}R\left(D+B\right)\right\} +{\rm Tr}\left\{ \tilde{W}^{\top}k\hat{W}\right\} 
\end{split}
\label{eq:equa49}
\end{equation}

\begin{equation}
\begin{split}\dot{V_{1}}+\dot{V_{2}}= & -c{\bf E}^{\top}\mathcal{Q}R{\bf E}-{\rm Tr}\left\{ \tilde{W}^{\top}\phi\left(x\right){\bf E}^{\top}\mathcal{M}RA\right\} +{\bf E}^{\top}\mathcal{M}\Delta \\
& +k{\rm Tr}\left\{ \tilde{W}^{\top}\hat{W}\right\} +{\bf E}^{\top}\mathcal{M}R\bigg(A\left(D+B\right)^{-1}R^{-1}\Phi_{2}\bar{\lambda}\\
& +\left(L+B\right)\left(\alpha-\underline{f}\left(x_{0},t\right)\right)\bigg)
\end{split}
\label{eq:equa50}
\end{equation}
\begin{equation}
\begin{split}\dot{V_{1}}+\dot{V_{2}}= & -c{\bf E}^{\top}\mathcal{Q}R{\bf E}-k{\rm Tr}\left\{ \tilde{W}^{\top}\tilde{W}\right\} -{\rm Tr}\left\{ \tilde{W}^{\top}\phi\left(x\right){\bf E}^{\top}\mathcal{M}RA\right\} \\
& +{\bf E}^{\top}\mathcal{M}\Delta +{\bf E}^{\top}\mathcal{M}R\bigg(A\left(D+B\right)^{-1}R^{-1}\Phi_{2}\bar{\lambda}\\
& +\left(L+B\right)\left(\alpha-\underline{f}\left(x_{0},t\right)\right)\bigg) +k{\rm Tr}\left\{ \tilde{W}^{\top}W\right\} 
\end{split}
\label{eq:equa51}
\end{equation}
Let $T_{M}=\alpha_{M}-\underline{f}_{M}$ such that
\begin{equation}
\begin{split}\dot{V_{1}}+\dot{V_{2}}\leq & -\left(c\underline{\sigma}\left(R\right)\underline{\sigma}\left(\mathcal{Q}\right)-\frac{\bar{\sigma}\left(\mathcal{M}\right)\bar{\sigma}\left(A\right)}{\underline{\sigma}\left(D+B\right)}\left\Vert \bar{\lambda}\right\Vert \right)\left\Vert {\bf E}\right\Vert ^{2}\\
& + \Bigg(-\bar{\sigma}\left(A\right)\phi_{M}\left(R\right)\left\Vert \tilde{W}\right\Vert _{F}+\frac{\bar{\sigma}\left(A\right)\left\Vert \Lambda\right\Vert {}_{F}\left\Vert \Phi_{1}\right\Vert }{\underline{\sigma}\left(D+B\right)}\\
& \hspace{10pt}+\bar{\sigma}\left(\Delta\right)+\bar{\sigma}\left(R\right)\bar{\sigma}\left(L+B\right)T_{M}\Bigg)\bar{\sigma}\left(\mathcal{M}\right)\left\Vert {\bf E}\right\Vert \\
& -k\left\Vert \tilde{W}\right\Vert _{F}^{2}+kW_{M}\left\Vert \tilde{W}\right\Vert _{F}
\end{split}
\label{eq:equa52}
\end{equation}
Now, the derivative of the third Lyapunov term $V_{3}$ is 
\begin{equation}
\begin{aligned}\dot{V_{3}}&=\frac{1}{2}{\rm Tr}\left\{ \dot{\Phi_{1}}M\Phi_{1}^{\top}+\Phi_{1}M\dot{\Phi_{1}}^{\top}\right\} \\
&={\rm Tr}\left\{ \Phi_{2}M\Phi_{1}^{\top}\right\} \end{aligned}
\label{eq:equa53}
\end{equation}
substituting \eqref{eq:equa24} in \eqref{eq:equa53} gives 
\begin{equation}
\begin{split}\dot{V_{3}}={\rm Tr}\left\{ \Phi_{1}\Lambda^{\top}M\Phi_{1}^{\top}\right\} +{\rm Tr}\left\{ {\bf E}l^{\top}M\Phi_{1}^{\top}\right\} \end{split}
\label{eq:equa54}
\end{equation}
which means
\begin{equation}
\begin{split}\dot{V_{3}} & ={\rm Tr}\left\{ \Phi_{1}(M\Lambda+\Lambda^{\top}M)\Phi_{1}^{\top}\right\} +{\rm Tr}\left\{ {\bf E}l^{\top}M\Phi_{1}^{\top}\right\} \\
& =-\frac{1}{2}\beta{\rm Tr}\left\{ \Phi_{1}\Phi_{1}^{\top}\right\} +{\rm Tr}\left\{ {\bf E}l^{\top}M\Phi_{1}^{\top}\right\} 
\end{split}
\label{eq:equa55}
\end{equation}
One can write the derivative of the complete form in \eqref{eq:equa45} as
\begin{equation}
\begin{split}\dot{V_{1}}+\dot{V_{2}}+\dot{V_{3}}\leq & -\left(c\underline{\sigma}\left(R\right)\underline{\sigma}\left(\mathcal{Q}\right)-\frac{\bar{\sigma}\left(\mathcal{M}\right)\bar{\sigma}\left(A\right)}{\underline{\sigma}\left(D+B\right)}\left\Vert \bar{\lambda}\right\Vert \right)\left\Vert {\bf E}\right\Vert ^{2}\\
& +\left(\bar{\sigma}\left(M\right)+\frac{\bar{\sigma}\left(\mathcal{M}\right)\bar{\sigma}\left(A\right)}{\underline{\sigma}\left(D+B\right)}\left\Vert \Lambda \right\Vert_{F}\right)\left\Vert \Phi_{1}\right\Vert \left\Vert {\bf E}\right\Vert\\
& +\bigg(\bar{\sigma}\left(R\right)\bar{\sigma}\left(L+B\right)T_{M}+\bar{\sigma}\left(\Delta\right)\\
& \hspace{10pt}-\phi_{M}\bar{\sigma}\left(R\right)\bar{\sigma}\left(A\right)\left\Vert \tilde{W}\right\Vert _{F}\bigg)\bar{\sigma}\left(\mathcal{M}\right)\left\Vert {\bf E}\right\Vert\\
& +kW_{M}\left\Vert \tilde{W}\right\Vert _{F} -k\left\Vert \tilde{W}\right\Vert _{F}^{2}-\frac{1}{2}\beta\left\Vert \Phi_{1}\right\Vert ^{2}
\end{split}
\label{eq:equa56}
\end{equation}
The result in \eqref{eq:equa56} can be simplified to
\small{
	\begin{equation}
	\begin{split}\dot{V}\leq & -\begin{bmatrix}\left\Vert \Phi_{1}\right\Vert  & \left\Vert \tilde{W}\right\Vert _{F} & \left\Vert {\bf E}\right\Vert \end{bmatrix}\begin{bmatrix}\frac{1}{2}\beta & 0 & g\\
	0 & k & \gamma\\
	g & \gamma & \mu
	\end{bmatrix}\begin{bmatrix}\left\Vert \Phi_{1}\right\Vert \\
	\left\Vert \tilde{W}\right\Vert _{F}\\
	\left\Vert {\bf E}\right\Vert 
	\end{bmatrix}\\
	& +\begin{bmatrix}0 & kW_{M} & \bar{\sigma}\left(\mathcal{M}\right)\left(\bar{\sigma}\left(R\right)\bar{\sigma}\left(L+B\right)T_{M}+\bar{\sigma}\left(\Delta\right)\right)\end{bmatrix}\begin{bmatrix}\left\Vert \Phi_{1}\right\Vert \\
	\left\Vert \tilde{W}\right\Vert _{F}\\
	\left\Vert {\bf E}\right\Vert 
	\end{bmatrix}
	\end{split}
	\label{eq:equa57}
	\end{equation}
}
\normalsize
with\\ $\gamma=-\frac{1}{2}\Phi\bar{\sigma}\left(\mathcal{M}\right)\bar{\sigma}\left(R\right)\bar{\sigma}\left(A\right)$,\\
$g=-\frac{1}{2}\left(\bar{\sigma}\left(M\right)+\frac{\bar{\sigma}\left(\mathcal{M}\right)\bar{\sigma}\left(A\right)}{\underline{\sigma}\left(D+B\right)}\left\Vert \Lambda \right\Vert_{F}\left\Vert \bar{\lambda}\right\Vert \right)$,\\ $\nu=\frac{\bar{\sigma}\left(\mathcal{M}\right)\bar{\sigma}\left(A\right)}{\underline{\sigma}\left(D+B\right)}\left\Vert \bar{\lambda}\right\Vert $, and\\ $\mu=\left(c\underline{\sigma}\left(R\right)\underline{\sigma}\left(\mathcal{Q}\right)-\frac{\bar{\sigma}\left(\mathcal{M}\right)\bar{\sigma}\left(A\right)}{\underline{\sigma}\left(D+B\right)}\right)$.\\
Let's write \eqref{eq:equa57} as 
\begin{equation}
\begin{split}\dot{V}\leq & -z^{\top} H z+h^{\top} z\end{split}
\label{eq:equa58}
\end{equation}
such that\\ $z=\begin{bmatrix}\left\Vert \Phi_{1}\right\Vert  & \left\Vert \tilde{W}\right\Vert _{F} & \left\Vert {\bf E}\right\Vert \end{bmatrix}^{\top}$,\\ $h=\begin{bmatrix}0 & kW_{M} & \bar{\sigma}\left(\mathcal{M}\right)\left(\bar{\sigma}\left(R\right)\bar{\sigma}\left(L+B\right)T_{M}+\bar{\sigma}\left(\Delta\right)\right)\end{bmatrix}^{\top}$,\\ $H = \begin{bmatrix}\frac{1}{2}\beta & 0 & g\\
0 & k & \gamma\\
g & \gamma & \mu
\end{bmatrix}$,\\ 
hence, $\dot{V}\leq0$ is only valid if $H$ is positive definite
there is 
\begin{equation}
\begin{split}\left\Vert z\right\Vert > & \frac{\left\Vert h\right\Vert }{\underline{\sigma}\left(H\right)}\end{split}
\label{eq:equa59}
\end{equation}
According to Sylvester's criterion, $H>0$ if 
\begin{enumerate}
	\item $\beta>0$ 
	\item $\beta k>0$ 
	\item $k\left(\beta\mu-2g^{2}\right)-\beta\gamma^{2}>0$ 
\end{enumerate}
Solving the foregoing equations proofs Theorem \ref{theom1}
\[
\begin{split}c>\frac{1}{\underline{\sigma}\left(\mathcal{Q}\right)\underline{\sigma}\left(R\right)}\left(\frac{\gamma^{2}}{k}+\frac{2}{\beta}g^{2}+\nu\right)\end{split}
\]
Then, assume 
\begin{equation}
\begin{split}\eta=\frac{kW_{M}+\bar{\sigma}\left(\mathcal{M}\right)\left(\bar{\sigma}\left(R\right)\bar{\sigma}\left(L+B\right)T_{M}+\bar{\sigma}\left(\Delta\right)\right)}{\underline{\sigma}\left(H\right)}\end{split}
\label{eq:equa60}
\end{equation}
We have $\dot{V}\leq0$ if $\left\Vert z\right\Vert >\eta$, according
to \eqref{eq:equa44}, we have
\small{
	\begin{equation}
	\begin{split} & \frac{1}{2}z^{\top}\begin{bmatrix}\underline{\sigma}\left(M\right) & 0 & 0\\
	0 & \bar{\sigma}\left(F\right) & 0\\
	0 & 0 & \underline{\sigma}\left(\mathcal{M}\right)
	\end{bmatrix}z\leq V\leq\frac{1}{2}z^{\top}\begin{bmatrix}\bar{\sigma}\left(M\right) & 0 & 0\\
	0 & \underline{\sigma}\left(F\right) & 0\\
	0 & 0 & \bar{\sigma}\left(\mathcal{M}\right)
	\end{bmatrix}z\end{split}
	\label{eq:equa61}
	\end{equation}
}
\normalsize
Define the appropriate variables matched with \eqref{eq:equa61} to
write 
\[
%\label{eq:eq45a}
\begin{split}\frac{1}{2}z^{\top} \underline{\mathcal{X}} z\leq V\leq\frac{1}{2}z^{\top} \bar{\mathcal{X}}z\end{split}
\]
Accordingly, it is equivalent to 
\begin{equation}
\begin{split}\frac{1}{2}\underline{\sigma}\left(\underline{\mathcal{X}}\right)\left\Vert z\right\Vert ^{2}\leq V\leq\frac{1}{2}\bar{\sigma}\left(\bar{\mathcal{X}}\right)\left\Vert z\right\Vert ^{2}\end{split}
\label{eq:equa62}
\end{equation}

\begin{equation}
\begin{split}V>\frac{1}{2}\bar{\sigma}\left(\bar{\mathcal{X}}\right)\frac{\left\Vert h\right\Vert ^{2}}{\underline{\sigma}^{2}\left(H\right)}\end{split}
\label{eq:equa63}
\end{equation}
Hence, based on Theorem 4.18 in \cite{khalil_nonlinear_2002}, for
any the initial value $z\left(t_{0}\right)$, there exists $T_{0}$
such that 
\begin{equation}
\begin{split}z\left(t\right)<\sqrt{\frac{\bar{\sigma}\left(\bar{\mathcal{X}}\right)}{\underline{\sigma}\left(\underline{\mathcal{X}}\right)}}\eta,\forall t\geq t_{0}+T_{0}\end{split}
\label{eq:equa64}
\end{equation}
The time $T_{0}$ can be evaluated by 
\begin{equation}
\begin{split}T_{0}=\frac{V\left(t_{0}\right)-\bar{\sigma}\left(\bar{\mathcal{X}}\right)\eta^{2}}{k}\end{split}
\label{eq:equa65}
\end{equation}
And according \eqref{eq:equa62}, one can find
\begin{equation}
\begin{split}\left\Vert z\right\Vert \leq\sqrt{\frac{2V}{\underline{\sigma}\left(\underline{\mathcal{X}}\right)}},\hspace{10pt}\left\Vert z\right\Vert \geq\sqrt{\frac{2V}{\bar{\sigma}\left(\bar{\mathcal{X}}\right)}}\end{split}
\label{eq:equa66}
\end{equation}
Therefore, $\dot{V}$ in \eqref{eq:equa57} can be given by
\begin{equation}
\dot{V}\leq-\tau_{1} V+\tau_{2} \sqrt{V}\label{eq:equa67}
\end{equation}
with $\tau_{1}=\frac{2\underline{\sigma}\left(H\right)}{\bar{\sigma}\left(\bar{\mathcal{X}}\right)}$
and $\tau_{2}=\frac{\sqrt{2}\left\Vert h\right\Vert }{\sqrt{\underline{\sigma}\left(\underline{\mathcal{X}}\right)}}$
yield 
\begin{equation}
\sqrt{V}\leq\sqrt{V\left(0\right)}+\frac{\tau_{2}}{\tau_{1}}\label{eq:equa68}
\end{equation}
It can be concluded that $\varepsilon$ is $\mathcal{L}_{\infty}$
for $t\geq t_{0}$ in a compact set $\Psi_{0}=\left\{\varepsilon\left(t\right)|\left\Vert \varepsilon\left(t\right)\right\Vert \leq r_{t_{0}}\right\}$.
Consequently, $e\left(t\right)$ will satisfy the prescribed performance
for all $t \geq 0$ if we start at $t=t_{0}$ within the prescribed functions.\\

\begin{rem}
	The control signal for high order MIMO systems can be written as
	{\small
	\begin{equation}
	\begin{aligned}
	\label{eq:equa69}
	u_{i} = & B_{i}^{-1}\left(-c{\bf E}_{i} - \hat{W_{i}}^{\top}\phi_{i}\left(x_{i}\right)\right) \\
	& - B_{i}^{-1}\left(\left(d_{i}+b_{i}\right)^{-1}\otimes \mathbb{I}_P\right)\Omega_{i}^{-1}\left( \lambda_{i}^{M_{p}-1}\varepsilon_{i}^{M_{p}} + \cdots + \lambda_{i}^{1}\varepsilon_{i}^{2}\right)
	\end{aligned}
	\end{equation} }
	and the estimated weights of NN can be defined by
	\begin{equation}
	\label{eq:equa70}
	\dot{\hat{W}}_{i} = F_{i}\phi_{i}{\bf E}_{i}^{\top}m_{i}\Omega_{i}\left(\left(d_{i}+b_{i}\right)\otimes \mathbb{I}_P\right) - kF_{i}\hat{W}_{i}
	\end{equation}
\end{rem}
In brief, the proposed algorithm of nonlinear high order agent dynamics as in equation \eqref{eq:equa1} can be summarized by
\begin{enumerate}
	\item[Step 1.] Select the setting parameters $\bar{\delta}_{i}^{p}$,
	$\underline{\delta}_{i}^{p}$, $\rho_{i,\infty}^{p}$, $\ell_{i}^{p}$, $\Pi_{i}$, $k$ and $c$.
	\item[Step 2.] Obtain the synchronized local error $e_{i}^{p}$ from \eqref{eq:equa5} or \eqref{eq:equa7}.
	\item[Step 3.] Evaluate the PPF $\rho_{i}^{p}$ from \eqref{eq:equa11}.
	\item[Step 4.] Evaluate $r_{i}^{p}$ from \eqref{eq:equa20}.
	\item[Step 5.] Obtain the transformed error from \eqref{eq:equa20} starting from $\varepsilon^{1}$ to $\varepsilon^{M_{p}}$.
	\item[Step 6.] Evaluate the metric error ${\bf E}_{i}$ from \eqref{eq:equa21} or \eqref{eq:equa22}.
	\item[Step 7.] The control signal $u_i$ from \eqref{eq:equa37} or \eqref{eq:equa69}.
	\item[Step 8.] Find $\hat{W}_{i}$ from \eqref{eq:equa39} or \eqref{eq:equa70}.
	\item[Step 9.] Go to Step 2.
\end{enumerate} 

%\end{proof}  

\section{Simulation Results} \label{Sec5}
\begin{figure}[h!]
	\centering
	\includegraphics[scale=0.5]{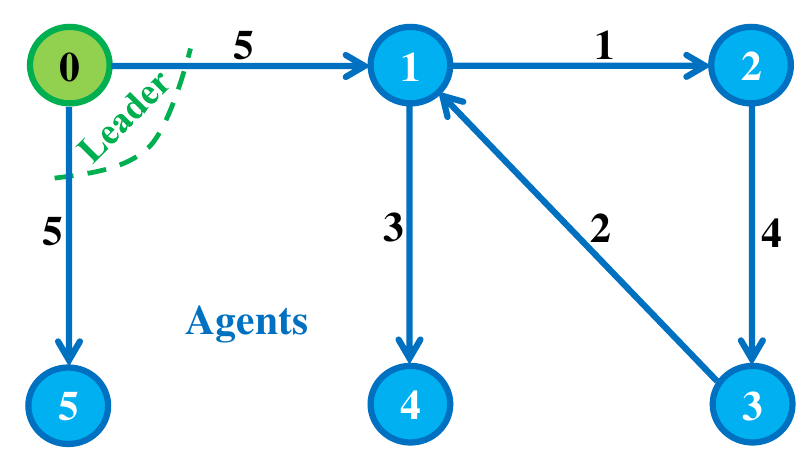}
	\caption{Five nodes with one leader in a connected graph.}
	\label{fig:fig2}
\end{figure}

{\bf Problem 1:} Consider the connected network in Figure \ref{fig:fig2} includes a group of five nodes with node (3) is connected to the leader node (0). The connected network is composed of 5 agents denoted by 1 to 5 with one leader denoted by 0 and the leader node is connected to agents 1 and 5 as in Figure \ref{fig:fig2}

\begin{figure*}[ht]
	\includegraphics[scale=0.45]{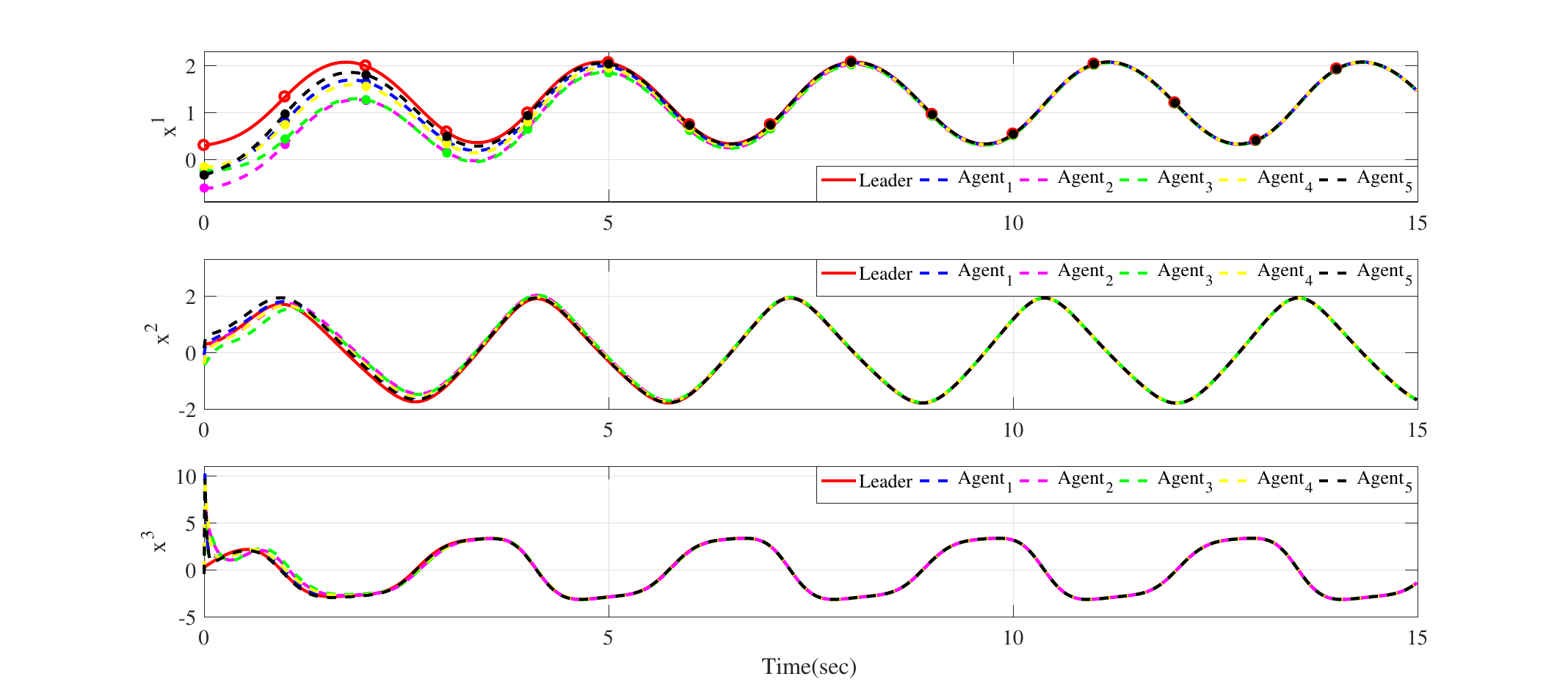}
	\caption{The output performance of high order nonlinear SISO networked system.}
	\label{fig:fig3}
\end{figure*}

\begin{figure*}[ht]
	\centering
	\includegraphics[scale=0.45]{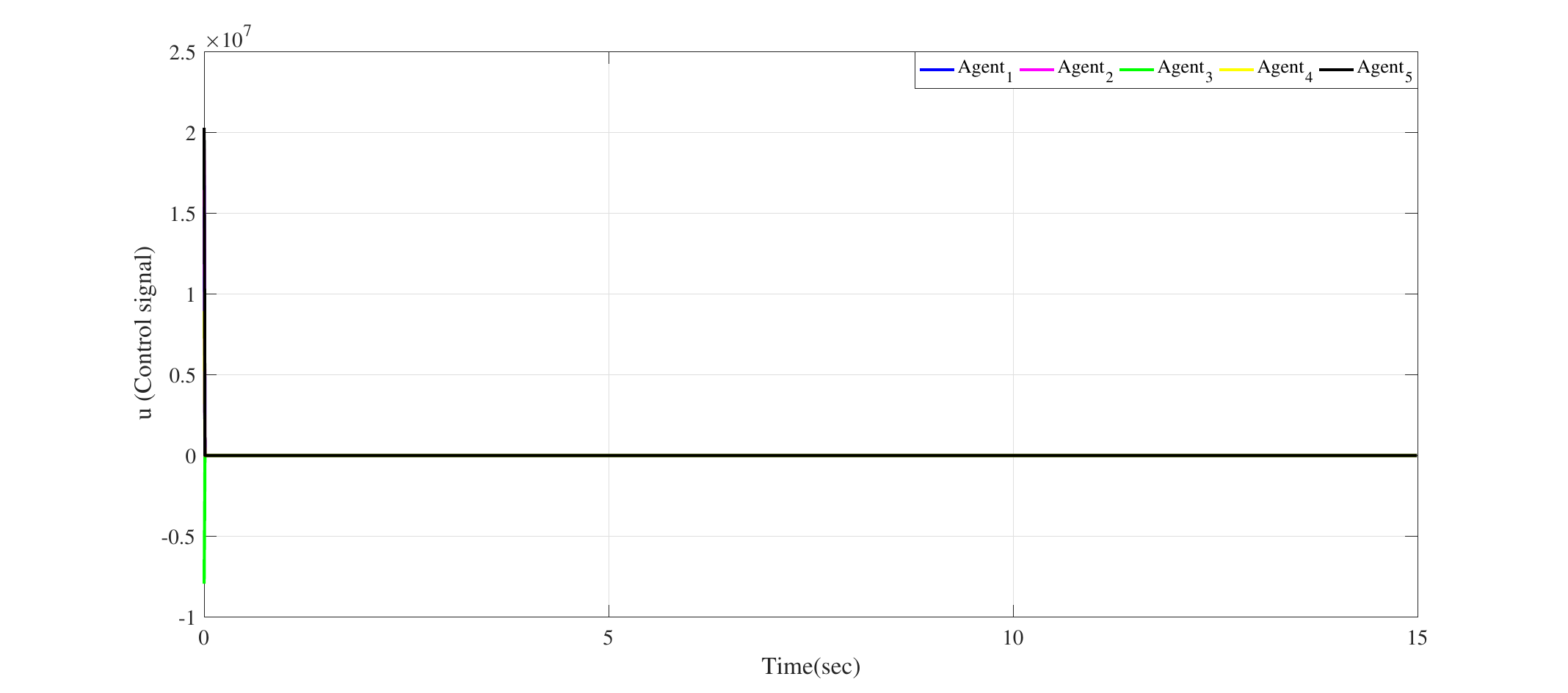}
	\caption{Control signal of high order nonlinear SISO networked system.}
	\label{fig:fig4}
\end{figure*}

\begin{figure*}[ht]
	\centering
	\includegraphics[scale=0.45]{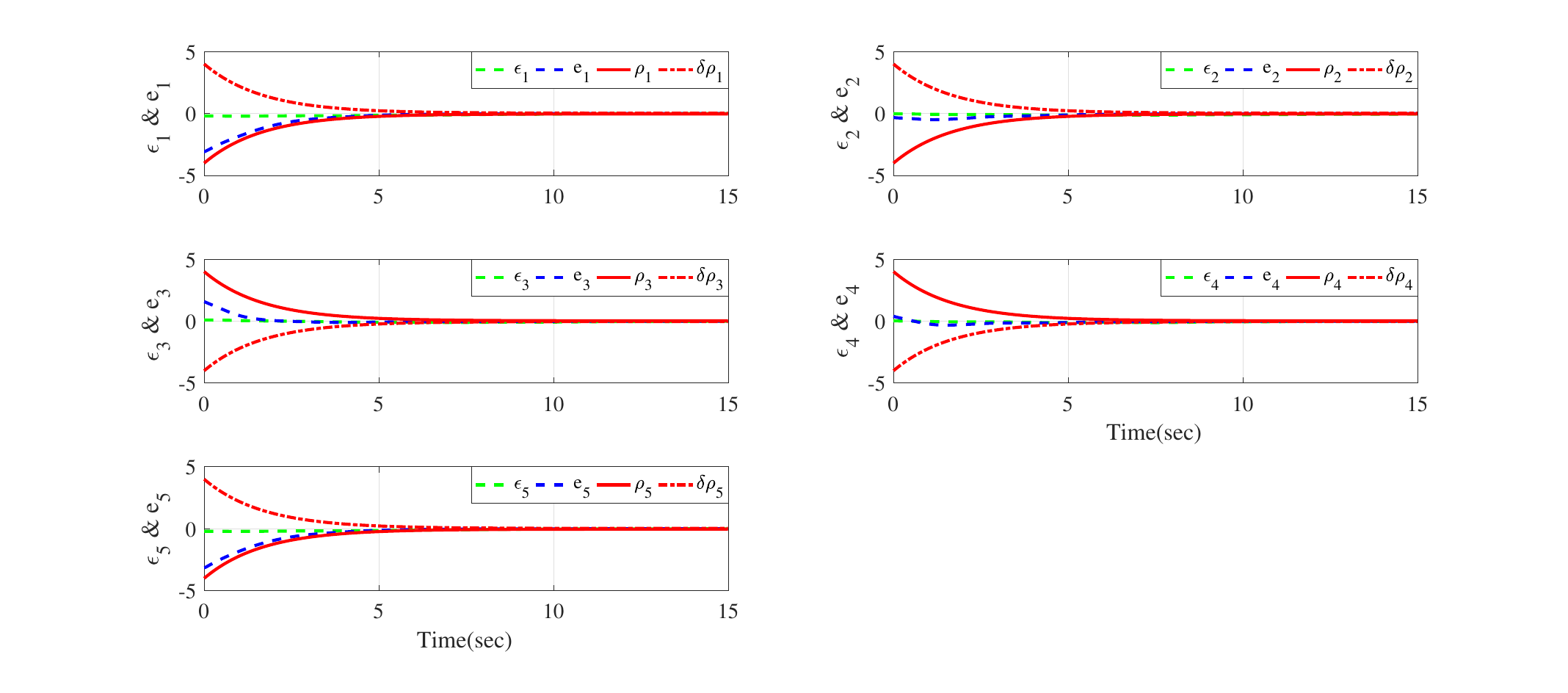}
	\caption{Error and transformed error of high order nonlinear SISO networked system.}
	\label{fig:fig5}
\end{figure*}

\begin{figure*}[ht]
	\centering
	\includegraphics[scale=0.5]{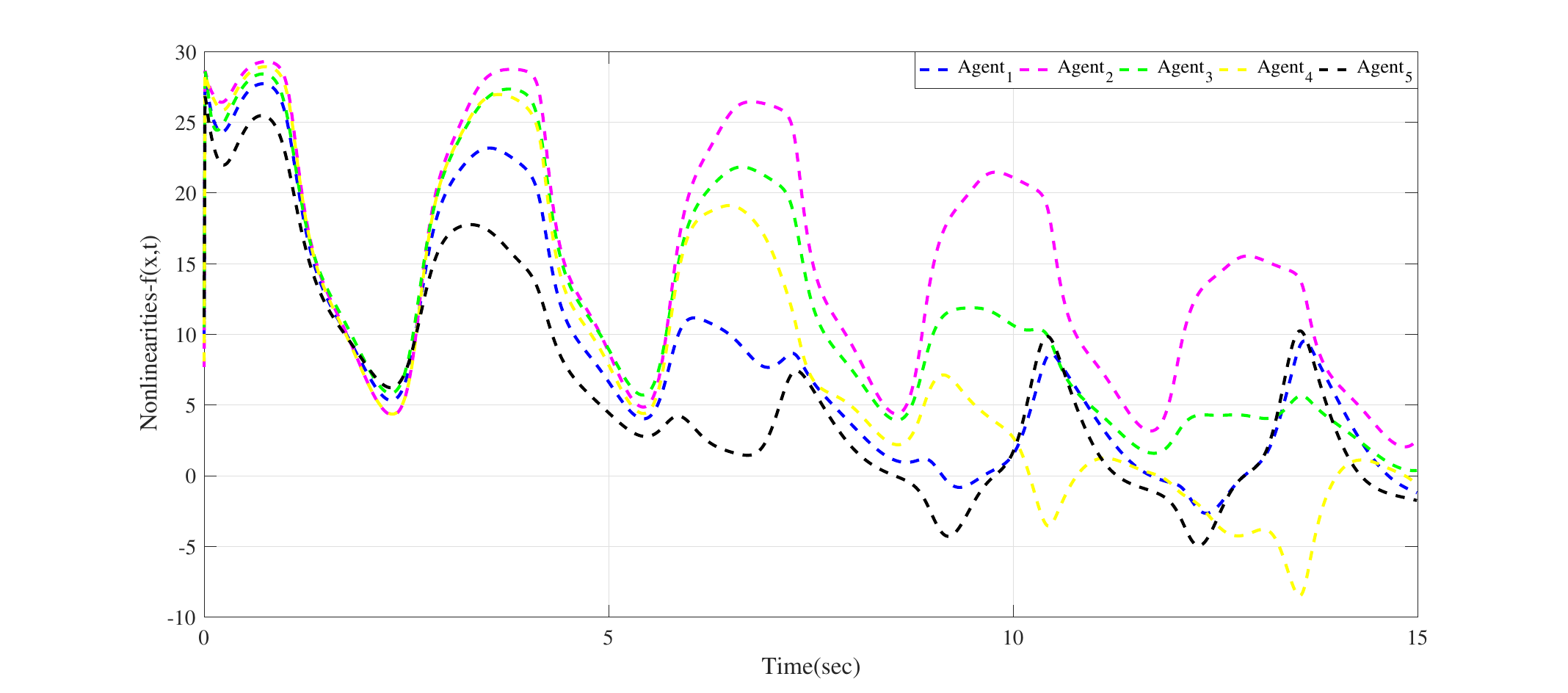}
	\caption{Nonlinearities compensation for each agent.}
	\label{fig:fig6}
\end{figure*}

The single-input single-output (SISO) dynamics of the graph are high order nonlinear such as
\begin{equation*}
\begin{aligned}
& \dot{x}_{i}^{1} = x_{i}^{2}\\
& \dot{x}_{i}^{2} = x_{i}^{3}\\
& \dot{x}_{i}^{3} = f_{i}\left(x_{i}\right) + u_{i}
\end{aligned}
\end{equation*}
such that $i = 1,2,\ldots,5$ with nonlinear dynamics
\begin{equation*}
\begin{aligned}
f_1 =& x_{1}^{2}{\rm sin}\left(x_{1}^{1}\right) + {\rm cos}\left(x_{1}^{3}\right)^2,\\
f_2 =& -\left(x_{2}^{1}\right)^2x_{2}^{2} + 0.01x_{2}^{1} - 0.01\left(x_{2}^{1}\right)^3,\\
f_3 =& x_{3}^{2} + {\rm sin}\left(x_{3}^{3}\right),\\
f_4 =& -3\left(x_{4}^{1}+x_{4}^{2}-1\right)^2\left(x_{4}^{1} + x_{4}^{2} + x_{4}^{3} - 1\right) - x_{4}^{3} + 0.5{\rm sin}\left(2t\right)\\
& +  {\rm cos}\left(2t\right),\\
f_5 =& {\rm cos}\left(x_{5}^{1}\right)
\end{aligned}
\end{equation*}
and the leader dynamics is 
\begin{equation*}
\begin{aligned}
\dot{x}_{0}^{1} = & x_{0}^{2}\\
\dot{x}_{0}^{2} = & x_{0}^{3}\\
\dot{x}_{0}^{3} = & -x_{0}^{2} - 2x_{0}^{3} + 1 + 3{\rm sin}\left(2t\right) + 6{\rm cos}\left(2t\right)\\
& - \frac{1}{3}\left(x_{0}^{1} + x_{0}^{2} - 1\right)\left(x_{0}^{1}+4x_{0}^{2} + 3x_{0}^{3}-1\right)
\end{aligned}
\end{equation*}
The setting parameters in this problem were selected as
$\rho_{\infty} = 0.03\times{\bf \underline{1}}_{5\times 1}$, $\rho_{0} = 4\times{\bf \underline{1}}_{5\times 1}$, $\ell = 0.6\times{\bf \underline{1}}_{5\times 1}$, $\Gamma = 0.05\mathbb{I}_{v_{i}\times v_{i}}$, $\bar{\delta} = 4\times{\bf \underline{1}}_{5\times 1}$, $\underline{\delta} = 4\times{\bf \underline{1}}_{5\times 1}$, $v_{i} = 6$, $c = 30\times{\bf \underline{1}}_{5\times 1}$, $k = 0.1$, 
$x_{0}\left(0\right) = [0.3,0.3,0.3]^{\top}$ is the initial vector of the nonlinear leader system, and the values of agents are\\ $x_{1}\left(0\right) = [-0.2850,-0.0821,-0.2126]^{\top}$,\\ $x_{2}\left(0\right) = [-0.6044,-0.3964,-0.0775]^{\top}$,\\ $x_{3}\left(0\right) = [-0.2110,-0.4237, -0.3253]^{\top}$,\\ $x_{4}\left(0\right) = [   -0.1501,-0.3986,-0.0050]^{\top}$,\\  $x_{5}\left(0\right) = [ -0.3281,0.1618,-0.4160]^{\top}$.\\      
Figure \ref{fig:fig3} presents the output performance of the proposed control and it shows smooth tracking performance for $x_1$ as well as for $x_2$ and $x_3$. It can be noticed that the networked system with unknown high nonlinear dynamics achieved consensus in the presence of high nonlinearities and time-varying disturbances. The control effort is depicted in Figure \ref{fig:fig4}. Figure \ref{fig:fig5} illustrates the systematic convergence of the synchronized error $e_i$ for $i=1,\ldots,5$ satisfying the predefined constraints and setting parameters imposed on the system. In fact, Figure \ref{fig:fig5} shows how the error started from a predefined large set and reduced systematically into the predefined small set prescribed by the value of $\rho_{\infty}$. Moreover, Figure \ref{fig:fig5} presents transformed error $\varepsilon_i$ associated to agent $i$. The nonlinear compensation of each agent is bounded and smooth as presented in Figure \ref{fig:fig6}. 

{\bf Problem 2:} Consider a second order dynamics of multi-input multi-output (MIMO) system with graph similar to Figure \ref{fig:fig2}. Each agent has 2 inputs and 2 outputs. The nonlinear dynamics are defined by
\begin{equation*}
\begin{aligned}
\begin{bmatrix} \ddot{x}_{i}^{1} \\ \ddot{x}_{i}^{2} \end{bmatrix}  & = \begin{bmatrix} f_{i}^{1}(x_{i},t) \\ f_{i}^{2}(x_{i},t) \end{bmatrix} + \psi_{i}\left(t\right)\begin{bmatrix}x_{i}^{1} \\ x_{i}^{2} \end{bmatrix} + \begin{bmatrix} D_{i}^{1}\left(t\right) \\ D_{i}^{2}\left(t\right) \end{bmatrix} + \begin{bmatrix} u_{i}^{1} \\ u_{i}^{2} \end{bmatrix}\\
y_{i,:} & = \begin{bmatrix} x_{i}^{1} & x_{i}^{2} \end{bmatrix}^{\top}
\end{aligned}
\end{equation*}
where        
\begin{equation*}
\begin{aligned}
f_{i}(x_{i}) = \begin{bmatrix}
a_{i}^{1}x_{i}^{2}x_{i}^{1}x_{i}^{1}\dot{x}_{i}^{2}+0.2{\rm sin}\left(a_{i}^{1}x_{i}^{1}\dot{x}_{i}^{1}\right)\\ -a_{i}^{2}x_{i}^{1}x_{i}^{2}\dot{x}_{i}^{1} -0.2a_{i}^{2}{\rm cos}\left(a_{i}^{2}x_{i}^{2}t\right)x_{i}^{1}\dot{x}_{i}^{2}
\end{bmatrix},
\end{aligned}
\end{equation*}

\begin{equation*}
\begin{aligned}
\psi_{i} = \begin{bmatrix}
3c_{i}^{1}{\rm sin}\left(0.5t\right) & 2c_{i}^{1}{\rm sin}\left(0.4c_{i}^{1}t\right){\rm cos}\left(0.3t\right)\\
0.9{\rm sin}\left(0.2c_{i}^{2}t\right) & 2.5{\rm sin}\left(0.3c_{i}^{2}t\right)+0.3{\rm cos}\left(t\right)
\end{bmatrix},
\end{aligned}
\end{equation*}

\begin{equation*}
\begin{aligned}
D_{i}\left(t\right) = \begin{bmatrix}
1+b_{i}^{1}{\rm sin}\left(b_{i}^{1}t\right)\\
1.2{\rm cos}\left(b_{i}^{2}t\right)
\end{bmatrix},
\end{aligned}
\end{equation*} 
and
\begin{equation*}
a =\begin{bmatrix}a_{i}^{1}\\a_{i}^{2}\end{bmatrix}= \begin{bmatrix}
1.5 & 0.5 & 0.7 & 1.3 & 0.7\\
0.5 & 1.4 & 0.1 & 1.3 & 2.4
\end{bmatrix}^{\top},
\end{equation*}
\begin{equation*}
b = \begin{bmatrix}
0.5 & 1.5 & 1.1 & 1.6 & 0.3\\
0.7 & 1.2 & 1.3 & 0.5 & 0.3
\end{bmatrix}^{\top},
\end{equation*}
\begin{equation*}
c = \begin{bmatrix}
1.5 & 2.5 & 0.5 & 1.7 & 0.7\\
0.5 & 1.7 & 1.1 & 0.3 & 0.4
\end{bmatrix}^{\top}
\end{equation*}

\begin{figure*}[ht]
	\centering
	\includegraphics[scale=0.45]{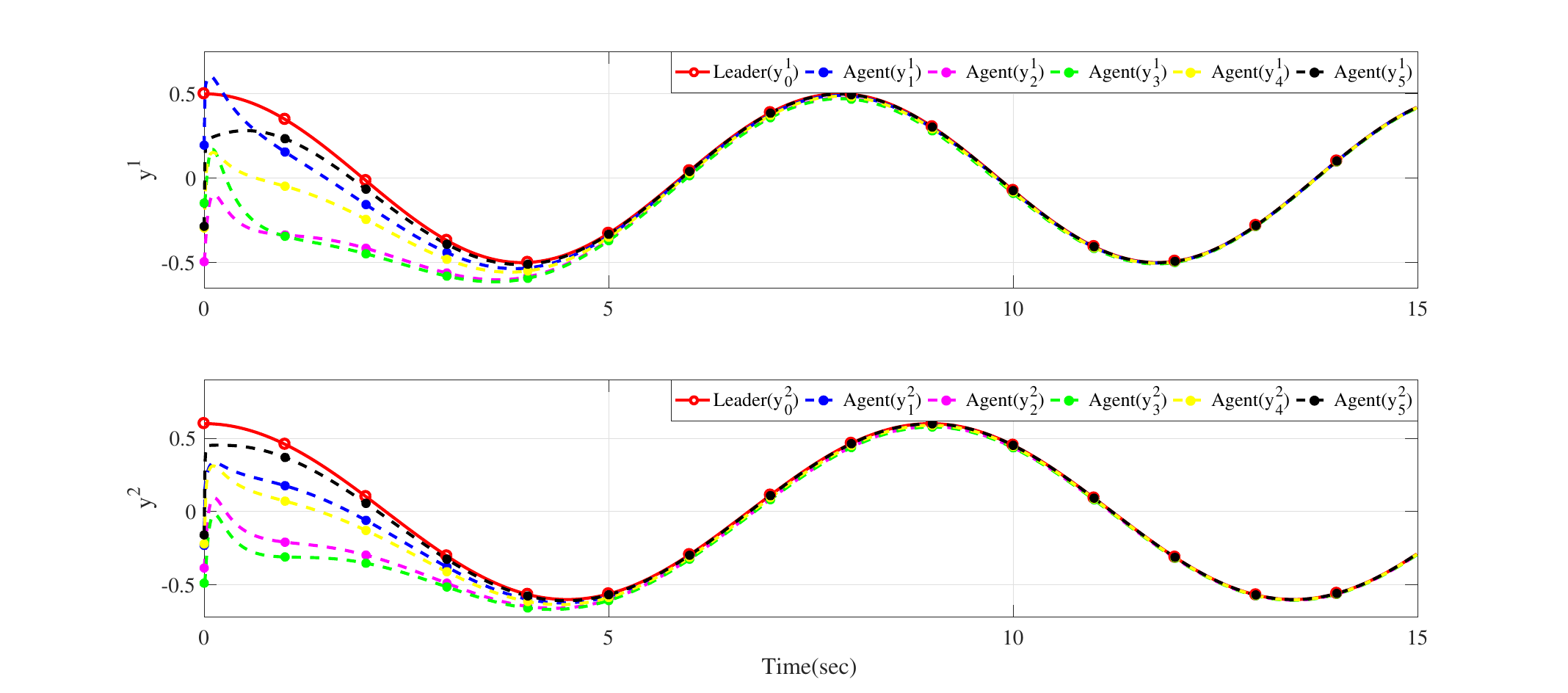}
	\caption{Output performance of MIMO nonlinear networked system for $y_1$ and $y_2$.}
	\label{fig:fig7}
\end{figure*}

\begin{figure*}[ht]
	\centering
	\includegraphics[scale=0.45]{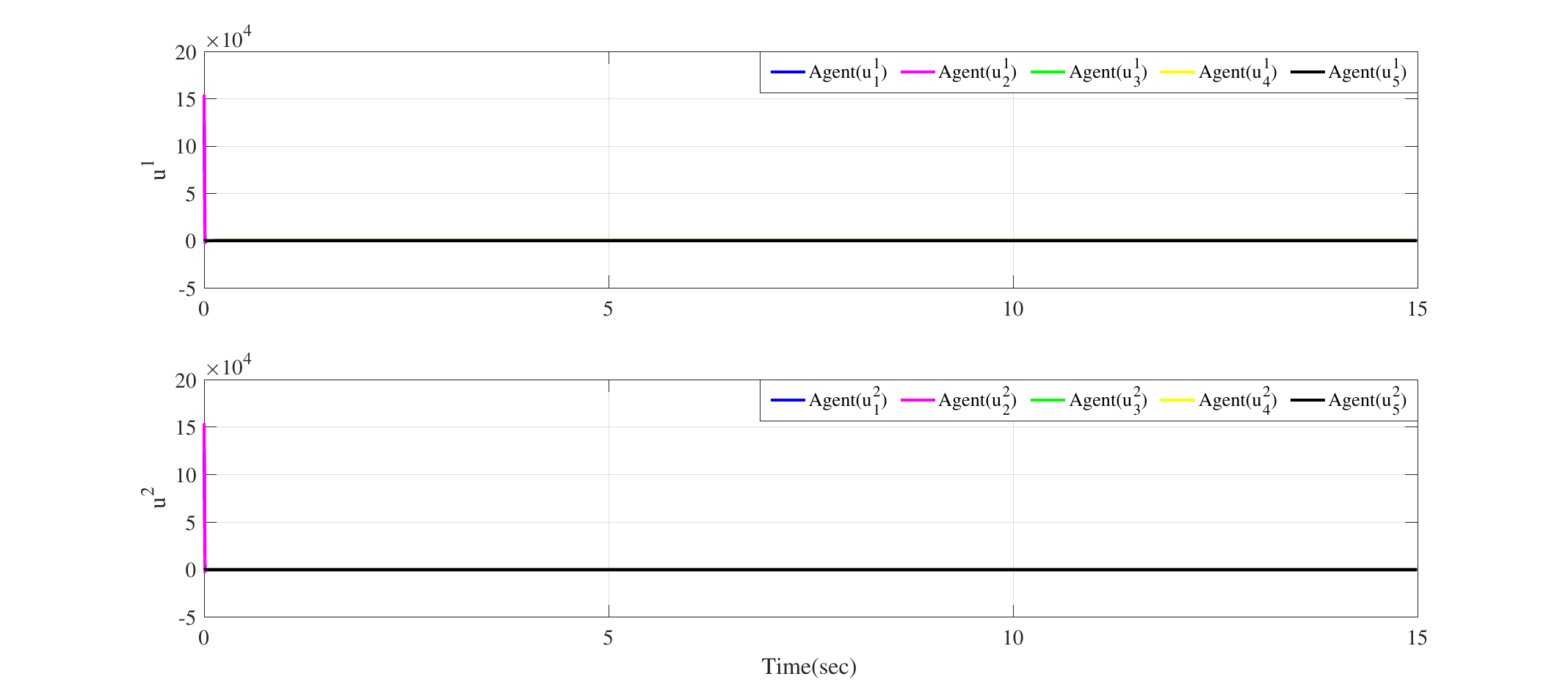}
	\caption{Control signal of MIMO nonlinear networked system for $u_1$ and $u_2$.}
	\label{fig:fig8}
\end{figure*}

\begin{figure*}[ht]
	\centering
	\includegraphics[scale=0.45]{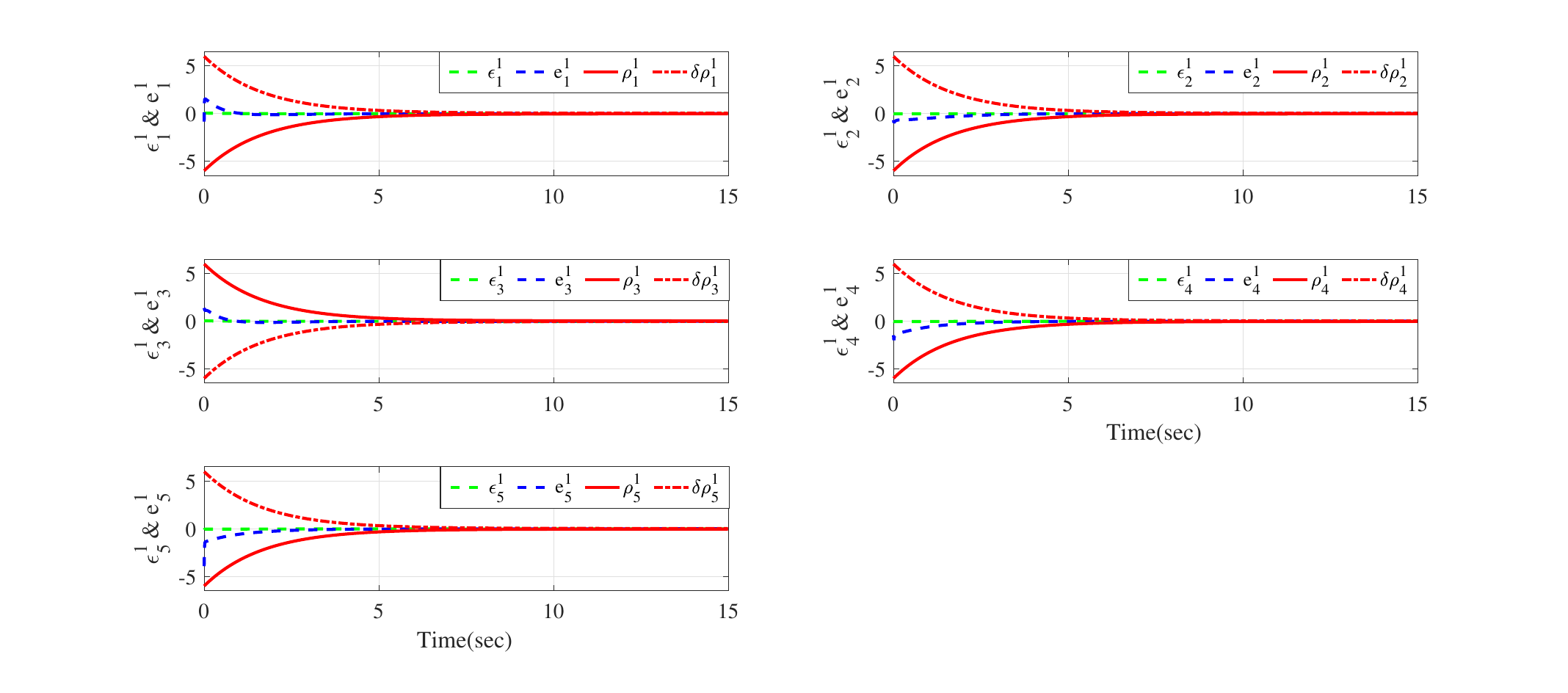}
	\caption{Error and transformed error of MIMO nonlinear networked system for the first output.} 
	\label{fig:fig9}    
\end{figure*}

\begin{figure*}[ht]
	\centering
	\includegraphics[scale=0.45]{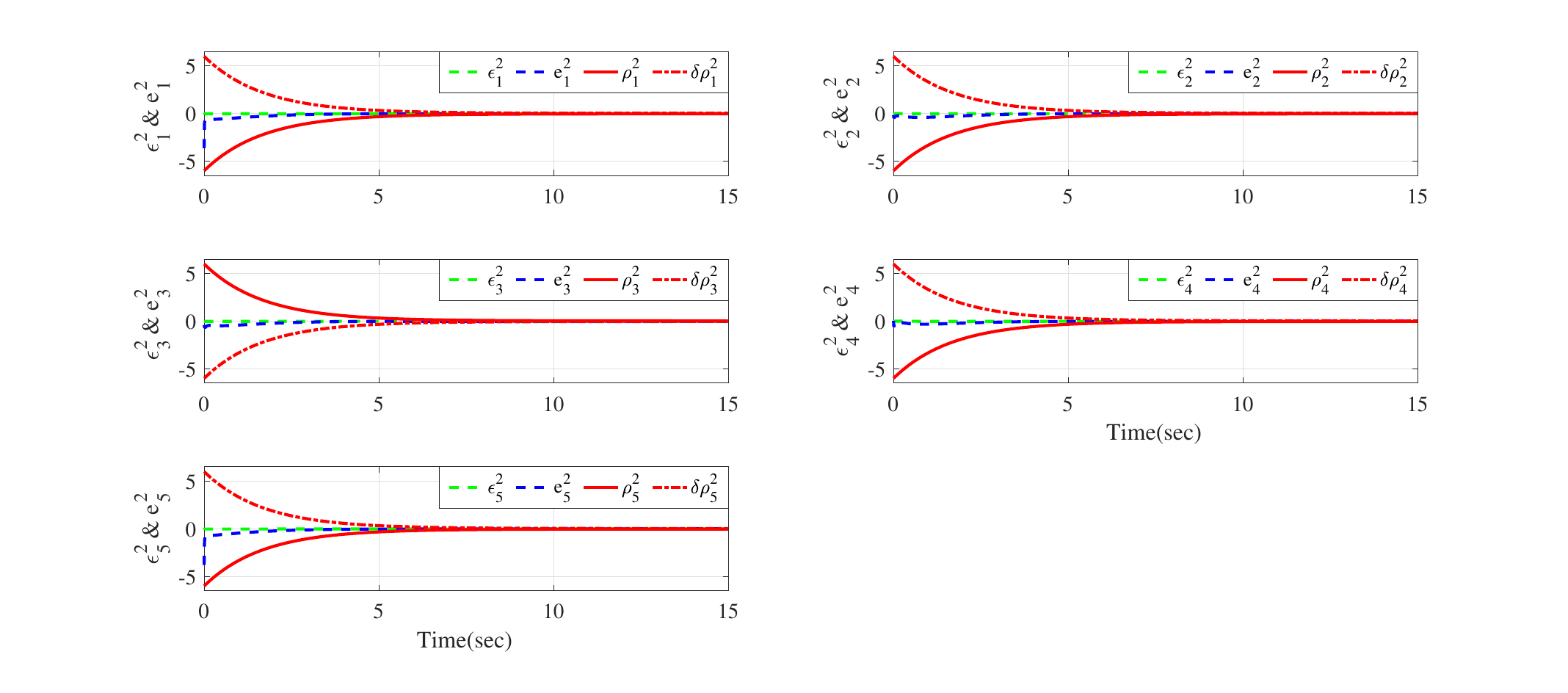}
	\caption{Error and transformed error of MIMO nonlinear networked system for the second output.}
	\label{fig:fig10}
\end{figure*}

\begin{figure*}[ht]
	\centering
	\includegraphics[scale=0.5]{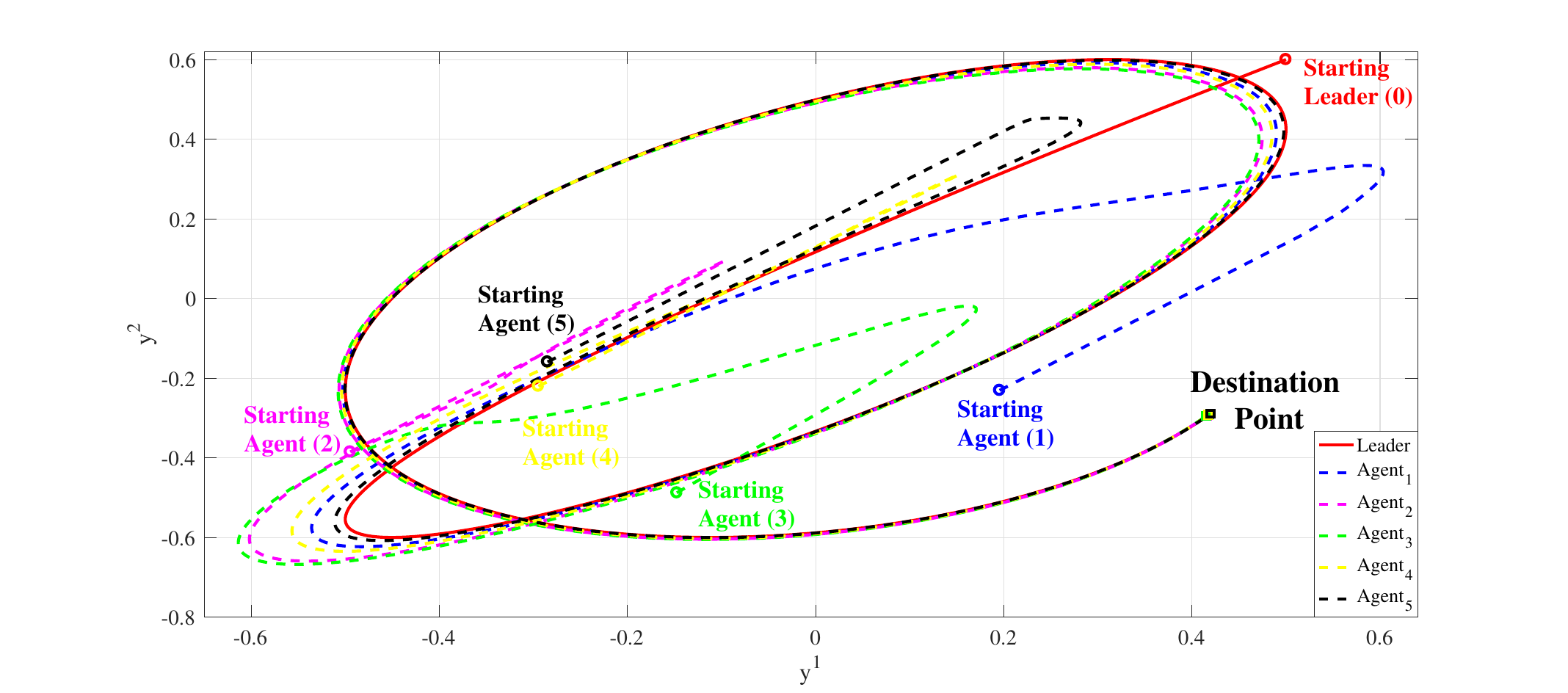}
	\caption{Phase plane plot for all the agents.}
	\label{fig:fig11}
\end{figure*}

The leader dynamics is $x_{0} = [0.5{\rm cos}\left(0.8t\right),0.6{\rm cos}\left(0.7t\right)]^{\top}$. In this problem, he setting parameters are $\rho_{0} = 6\times{\bf {\bf \underline{1}}_{5\times 2}}$, $\rho_{\infty} = 0.03\times{\bf {\bf \underline{1}}_{5\times 2}}$, $\ell = 0.6\times{\bf {\bf \underline{1}}_{5\times 2}}$, $\Gamma = 0.05\mathbb{I}_{5 \times 2}$, $\bar{\delta} = 6\times{\bf {\bf \underline{1}}_{5\times 2}}$, $\underline{\delta} = 6\times{\bf {\bf \underline{1}}_{5\times 2}}$, $c = 300\mathbb{I}_{5 \times 2}$, $k = 0.1\mathbb{I}_{5 \times 2}$. Initial conditions of\\ $x_1\left(0\right) = \left[0.1956,-0.2307\right]^{\top}$,\\ $x_2\left(0\right) = \left[-0.4947, -0.3852\right]^{\top}$,\\  $x_3\left(0\right) = \left[-0.1475, -0.4880\right]^{\top}$,\\ $x_4\left(0\right) = \left[-0.2947, -0.2203\right]^{\top}$,\\ $x_5\left(0\right) = \left[ -0.2850, -0.1593\right]^{\top}$\\ and $\dot{x}_{i}\left(0\right)=\left[0,0\right]^{\top}, i=1\ldots,5$.  Finally, number of neurons is $v_{i} = 50$.\\
The effectiveness of the proposed control algorithm against unknown nonlinearities and time variant components is examined in this problem. The output performance illustrates the robustness and high tracking capabilities of control algorithm. The output performance for high order MIMO case is given in Figure \ref{fig:fig7}. The output performance in Figure \ref{fig:fig7} moved from random initialization to consensus with the leader softly. The control input of system dynamics in the connected graph is shown in Figure \ref{fig:fig8}. Each of transformed and tracking errors are illustrated in figures \ref{fig:fig9} and \ref{fig:fig10}. One can clearly observe from these figures the quality of the control approach and its efficiency in ensuring synchronization with prescribed performance characteristics. Clearly, errors obeyed the predefined constraints and settings and the error started within a predefined large set to end within a predefined small set. Finally, the phase plane is presented in Figure \ref{fig:fig11} to show the random initialization of each agent with systematic consensus up to out destination or desired point.


\section{Conclusion} \label{Sec6}
Neuro-adaptive distributed control with prescribed performance of higher order nonlinear affine multi-agent systems with full-state synchronization has been proposed. Neural network is employed to estimate unknown nonlinear  dynamics of each node. The control signal has been chosen to both respect the digraph and ensure stability. The Proposed controller successfully allowed the nodes to synchronize the leader trajectory and satisfy at any point in time the desired performance with small residual errors. Also, the controller guarantees tracking the leaders' states with a synchronization error within a predefined time varying constraints. Prescribed performance controller was designed based on robust neuro-adaptive approach. Lyapunov-based stability proofs establish that the synchronization error of each node is UUB. Simulation examples related to single-input single-output and multi-input multi-output consider high-order dynamics with unknown nonlinearities and time variant components. However, the class of systems addressed in this paper possesses an input function with constant values. More realistic systems with general nonlinear input functions can be considered in future work. In addition, hard actuator nonlinearities and especially saturation input functions could be studied in future research.

\section*{Acknowledgment}
The authors received no financial support for the research, authorship, and/or publication of this article.
%\ifCLASSOPTIONcaptionsoff %\bibliographystyle{spmpsci}      % mathematics and physical sciences
\bibliographystyle{apacite}
\bibliography{bib_Nuero_Hashim}

\begin{thebibliography}{}

\bibitem [\protect \citeauthoryear {%
Bechlioulis%
\ \BBA {} Rovithakis%
}{%
Bechlioulis%
\ \BBA {} Rovithakis%
}{%
{\protect \APACyear {2008}}%
}]{%
bechlioulis_robust_2008}
\APACinsertmetastar {%
bechlioulis_robust_2008}%
\begin{APACrefauthors}%
Bechlioulis, C\BPBI P.%
\BCBT {}\ \BBA {} Rovithakis, G\BPBI A.%
\end{APACrefauthors}%
\unskip\
\newblock
\APACrefYearMonthDay{2008}{}{}.
\newblock
{\BBOQ}\APACrefatitle {Robust adaptive control of feedback linearizable {MIMO}
  nonlinear systems with prescribed performance} {Robust adaptive control of
  feedback linearizable {MIMO} nonlinear systems with prescribed
  performance}.{\BBCQ}
\newblock
\APACjournalVolNumPages{{IEEE} Transactions on Automatic
  Control}{53}{9}{2090--2099}.
\PrintBackRefs{\CurrentBib}

\bibitem [\protect \citeauthoryear {%
Bu%
, Wu%
, Huang%
\BCBL {}\ \BBA {} Wei%
}{%
Bu%
\ \protect \BOthers {.}}{%
{\protect \APACyear {2016}}%
}]{%
bu2016robust}
\APACinsertmetastar {%
bu2016robust}%
\begin{APACrefauthors}%
Bu, X.%
, Wu, X.%
, Huang, J.%
\BCBL {}\ \BBA {} Wei, D.%
\end{APACrefauthors}%
\unskip\
\newblock
\APACrefYearMonthDay{2016}{}{}.
\newblock
{\BBOQ}\APACrefatitle {Robust estimation-free prescribed performance
  back-stepping control of air-breathing hypersonic vehicles without affine
  models} {Robust estimation-free prescribed performance back-stepping control
  of air-breathing hypersonic vehicles without affine models}.{\BBCQ}
\newblock
\APACjournalVolNumPages{International Journal of Control}{89}{11}{2185--2200}.
\PrintBackRefs{\CurrentBib}

\bibitem [\protect \citeauthoryear {%
Cao%
\ \BBA {} Ren%
}{%
Cao%
\ \BBA {} Ren%
}{%
{\protect \APACyear {2012}}%
}]{%
cao_distributed_2012}
\APACinsertmetastar {%
cao_distributed_2012}%
\begin{APACrefauthors}%
Cao, Y.%
\BCBT {}\ \BBA {} Ren, W.%
\end{APACrefauthors}%
\unskip\
\newblock
\APACrefYearMonthDay{2012}{}{}.
\newblock
{\BBOQ}\APACrefatitle {Distributed coordinated tracking with reduced
  interaction via a variable structure approach} {Distributed coordinated
  tracking with reduced interaction via a variable structure approach}.{\BBCQ}
\newblock
\APACjournalVolNumPages{IEEE Transactions on Automatic Control}{57}{1}{33--48}.
\PrintBackRefs{\CurrentBib}

\bibitem [\protect \citeauthoryear {%
Chopra%
\ \BBA {} Spong%
}{%
Chopra%
\ \BBA {} Spong%
}{%
{\protect \APACyear {2006}}%
}]{%
chopra_passivity-based_2006}
\APACinsertmetastar {%
chopra_passivity-based_2006}%
\begin{APACrefauthors}%
Chopra, N.%
\BCBT {}\ \BBA {} Spong, M\BPBI W.%
\end{APACrefauthors}%
\unskip\
\newblock
\APACrefYearMonthDay{2006}{}{}.
\newblock
{\BBOQ}\APACrefatitle {Passivity-based control of multi-agent systems}
  {Passivity-based control of multi-agent systems}.{\BBCQ}
\newblock
\BIn{} \APACrefbtitle {Advances in robot control} {Advances in robot control}\
  (\BPGS\ 107--134).
\newblock
\APACaddressPublisher{}{Springer}.
\PrintBackRefs{\CurrentBib}

\bibitem [\protect \citeauthoryear {%
Cotter%
}{%
Cotter%
}{%
{\protect \APACyear {1989}}%
}]{%
cotter_stone-weierstrass_1989}
\APACinsertmetastar {%
cotter_stone-weierstrass_1989}%
\begin{APACrefauthors}%
Cotter, N\BPBI E.%
\end{APACrefauthors}%
\unskip\
\newblock
\APACrefYearMonthDay{1989}{}{}.
\newblock
{\BBOQ}\APACrefatitle {The {Stone}-{Weierstrass} theorem and its application to
  neural networks.} {The {Stone}-{Weierstrass} theorem and its application to
  neural networks.}{\BBCQ}
\newblock
\APACjournalVolNumPages{IEEE transactions on neural networks/a publication of
  the IEEE Neural Networks Council}{1}{4}{290--295}.
\PrintBackRefs{\CurrentBib}

\bibitem [\protect \citeauthoryear {%
Das%
\ \BBA {} Lewis%
}{%
Das%
\ \BBA {} Lewis%
}{%
{\protect \APACyear {2010}}%
}]{%
das_distributed_2010}
\APACinsertmetastar {%
das_distributed_2010}%
\begin{APACrefauthors}%
Das, A.%
\BCBT {}\ \BBA {} Lewis, F\BPBI L.%
\end{APACrefauthors}%
\unskip\
\newblock
\APACrefYearMonthDay{2010}{}{}.
\newblock
{\BBOQ}\APACrefatitle {Distributed adaptive control for synchronization of
  unknown nonlinear networked systems} {Distributed adaptive control for
  synchronization of unknown nonlinear networked systems}.{\BBCQ}
\newblock
\APACjournalVolNumPages{Automatica}{46}{12}{2014--2021}.
\PrintBackRefs{\CurrentBib}

\bibitem [\protect \citeauthoryear {%
El-Ferik%
, Hashim%
\BCBL {}\ \BBA {} Lewis%
}{%
El-Ferik%
\ \protect \BOthers {.}}{%
{\protect \APACyear {2017}}%
}]{%
el2017neuro}
\APACinsertmetastar {%
el2017neuro}%
\begin{APACrefauthors}%
El-Ferik, S.%
, Hashim, H\BPBI A.%
\BCBL {}\ \BBA {} Lewis, F\BPBI L.%
\end{APACrefauthors}%
\unskip\
\newblock
\APACrefYearMonthDay{2017}{}{}.
\newblock
{\BBOQ}\APACrefatitle {Neuro-Adaptive Distributed Control With Prescribed
  Performance for the Synchronization of Unknown Nonlinear Networked Systems}
  {Neuro-adaptive distributed control with prescribed performance for the
  synchronization of unknown nonlinear networked systems}.{\BBCQ}
\newblock
\APACjournalVolNumPages{IEEE Transactions on Systems, Man, and Cybernetics:
  Systems}{}{}{}.
\PrintBackRefs{\CurrentBib}

\bibitem [\protect \citeauthoryear {%
El-Ferik%
, Qureshi%
\BCBL {}\ \BBA {} Lewis%
}{%
El-Ferik%
\ \protect \BOthers {.}}{%
{\protect \APACyear {2014}}%
}]{%
el-ferik_neuro-adaptive_2014}
\APACinsertmetastar {%
el-ferik_neuro-adaptive_2014}%
\begin{APACrefauthors}%
El-Ferik, S.%
, Qureshi, A.%
\BCBL {}\ \BBA {} Lewis, F\BPBI L.%
\end{APACrefauthors}%
\unskip\
\newblock
\APACrefYearMonthDay{2014}{}{}.
\newblock
{\BBOQ}\APACrefatitle {Neuro-adaptive cooperative tracking control of unknown
  higher-order affine nonlinear systems} {Neuro-adaptive cooperative tracking
  control of unknown higher-order affine nonlinear systems}.{\BBCQ}
\newblock
\APACjournalVolNumPages{Automatica}{50}{3}{798--808}.
\PrintBackRefs{\CurrentBib}

\bibitem [\protect \citeauthoryear {%
Fax%
\ \BBA {} Murray%
}{%
Fax%
\ \BBA {} Murray%
}{%
{\protect \APACyear {2004}}%
}]{%
fax_information_2004}
\APACinsertmetastar {%
fax_information_2004}%
\begin{APACrefauthors}%
Fax, J\BPBI A.%
\BCBT {}\ \BBA {} Murray, R\BPBI M.%
\end{APACrefauthors}%
\unskip\
\newblock
\APACrefYearMonthDay{2004}{}{}.
\newblock
{\BBOQ}\APACrefatitle {Information flow and cooperative control of vehicle
  formations} {Information flow and cooperative control of vehicle
  formations}.{\BBCQ}
\newblock
\APACjournalVolNumPages{IEEE Transactions on Automatic
  Control}{49}{9}{1465--1476}.
\PrintBackRefs{\CurrentBib}

\bibitem [\protect \citeauthoryear {%
Hashim%
, El-Ferik%
\BCBL {}\ \BBA {} Lewis%
}{%
Hashim%
\ \protect \BOthers {.}}{%
{\protect \APACyear {2017}}%
}]{%
Hashim2017adaptive}
\APACinsertmetastar {%
Hashim2017adaptive}%
\begin{APACrefauthors}%
Hashim, H\BPBI A.%
, El-Ferik, S.%
\BCBL {}\ \BBA {} Lewis, F\BPBI L.%
\end{APACrefauthors}%
\unskip\
\newblock
\APACrefYearMonthDay{2017}{}{}.
\newblock
{\BBOQ}\APACrefatitle {Adaptive synchronisation of unknown nonlinear networked
  systems with prescribed performance} {Adaptive synchronisation of unknown
  nonlinear networked systems with prescribed performance}.{\BBCQ}
\newblock
\APACjournalVolNumPages{International Journal of Systems
  Science}{48}{4}{885--898}.
\PrintBackRefs{\CurrentBib}

\bibitem [\protect \citeauthoryear {%
Hornik%
, Stinchcombe%
\BCBL {}\ \BBA {} White%
}{%
Hornik%
\ \protect \BOthers {.}}{%
{\protect \APACyear {1989}}%
}]{%
hornik_multilayer_1989}
\APACinsertmetastar {%
hornik_multilayer_1989}%
\begin{APACrefauthors}%
Hornik, K.%
, Stinchcombe, M.%
\BCBL {}\ \BBA {} White, H.%
\end{APACrefauthors}%
\unskip\
\newblock
\APACrefYearMonthDay{1989}{}{}.
\newblock
{\BBOQ}\APACrefatitle {Multilayer feedforward networks are universal
  approximators} {Multilayer feedforward networks are universal
  approximators}.{\BBCQ}
\newblock
\APACjournalVolNumPages{Neural networks}{2}{5}{359--366}.
\PrintBackRefs{\CurrentBib}

\bibitem [\protect \citeauthoryear {%
Khalil%
}{%
Khalil%
}{%
{\protect \APACyear {2002}}%
}]{%
khalil_nonlinear_2002}
\APACinsertmetastar {%
khalil_nonlinear_2002}%
\begin{APACrefauthors}%
Khalil, H\BPBI K.%
\end{APACrefauthors}%
\unskip\
\newblock
\APACrefYearMonthDay{2002}{}{}.
\newblock
{\BBOQ}\APACrefatitle {Nonlinear systems, 3rd} {Nonlinear systems, 3rd}.{\BBCQ}
\newblock
\APACjournalVolNumPages{New Jewsey, Prentice Hall}{9}{}{}.
\PrintBackRefs{\CurrentBib}

\bibitem [\protect \citeauthoryear {%
Khoo%
, Xie%
\BCBL {}\ \BBA {} Man%
}{%
Khoo%
\ \protect \BOthers {.}}{%
{\protect \APACyear {2009}}%
}]{%
khoo_robust_2009}
\APACinsertmetastar {%
khoo_robust_2009}%
\begin{APACrefauthors}%
Khoo, S.%
, Xie, L.%
\BCBL {}\ \BBA {} Man, Z.%
\end{APACrefauthors}%
\unskip\
\newblock
\APACrefYearMonthDay{2009}{}{}.
\newblock
{\BBOQ}\APACrefatitle {Robust finite-time consensus tracking algorithm for
  multirobot systems} {Robust finite-time consensus tracking algorithm for
  multirobot systems}.{\BBCQ}
\newblock
\APACjournalVolNumPages{Mechatronics, IEEE/ASME Transactions
  on}{14}{2}{219--228}.
\PrintBackRefs{\CurrentBib}

\bibitem [\protect \citeauthoryear {%
F\BPBI L.~Lewis%
, Zhang%
, Hengster-Movric%
\BCBL {}\ \BBA {} Das%
}{%
F\BPBI L.~Lewis%
\ \protect \BOthers {.}}{%
{\protect \APACyear {2013}}%
}]{%
lewis_cooperative_2013}
\APACinsertmetastar {%
lewis_cooperative_2013}%
\begin{APACrefauthors}%
Lewis, F\BPBI L.%
, Zhang, H.%
, Hengster-Movric, K.%
\BCBL {}\ \BBA {} Das, A.%
\end{APACrefauthors}%
\unskip\
\newblock
\APACrefYear{2013}.
\newblock
\APACrefbtitle {Cooperative control of multi-agent systems: optimal and
  adaptive design approaches} {Cooperative control of multi-agent systems:
  optimal and adaptive design approaches}.
\newblock
\APACaddressPublisher{}{Springer Science \& Business Media}.
\PrintBackRefs{\CurrentBib}

\bibitem [\protect \citeauthoryear {%
F\BPBI W.~Lewis%
, Jagannathan%
\BCBL {}\ \BBA {} Yesildirak%
}{%
F\BPBI W.~Lewis%
\ \protect \BOthers {.}}{%
{\protect \APACyear {1998}}%
}]{%
lewis_neural_1998}
\APACinsertmetastar {%
lewis_neural_1998}%
\begin{APACrefauthors}%
Lewis, F\BPBI W.%
, Jagannathan, S.%
\BCBL {}\ \BBA {} Yesildirak, A.%
\end{APACrefauthors}%
\unskip\
\newblock
\APACrefYear{1998}.
\newblock
\APACrefbtitle {Neural network control of robot manipulators and non-linear
  systems} {Neural network control of robot manipulators and non-linear
  systems}.
\newblock
\APACaddressPublisher{}{CRC Press}.
\PrintBackRefs{\CurrentBib}

\bibitem [\protect \citeauthoryear {%
Li%
, Wang%
\BCBL {}\ \BBA {} Chen%
}{%
Li%
\ \protect \BOthers {.}}{%
{\protect \APACyear {2004}}%
}]{%
li_pinning_2004}
\APACinsertmetastar {%
li_pinning_2004}%
\begin{APACrefauthors}%
Li, X.%
, Wang, X.%
\BCBL {}\ \BBA {} Chen, G.%
\end{APACrefauthors}%
\unskip\
\newblock
\APACrefYearMonthDay{2004}{}{}.
\newblock
{\BBOQ}\APACrefatitle {Pinning a complex dynamical network to its equilibrium}
  {Pinning a complex dynamical network to its equilibrium}.{\BBCQ}
\newblock
\APACjournalVolNumPages{Circuits and Systems I: Regular Papers, IEEE
  Transactions on}{51}{10}{2074--2087}.
\PrintBackRefs{\CurrentBib}

\bibitem [\protect \citeauthoryear {%
Liao%
, Lu%
\BCBL {}\ \BBA {} Liu%
}{%
Liao%
\ \protect \BOthers {.}}{%
{\protect \APACyear {2016}}%
}]{%
liao2016cooperative}
\APACinsertmetastar {%
liao2016cooperative}%
\begin{APACrefauthors}%
Liao, F.%
, Lu, Y.%
\BCBL {}\ \BBA {} Liu, H.%
\end{APACrefauthors}%
\unskip\
\newblock
\APACrefYearMonthDay{2016}{}{}.
\newblock
{\BBOQ}\APACrefatitle {Cooperative optimal preview tracking control of
  continuous-time multi-agent systems} {Cooperative optimal preview tracking
  control of continuous-time multi-agent systems}.{\BBCQ}
\newblock
\APACjournalVolNumPages{International Journal of Control}{89}{10}{2019--2028}.
\PrintBackRefs{\CurrentBib}

\bibitem [\protect \citeauthoryear {%
Mohamed%
}{%
Mohamed%
}{%
{\protect \APACyear {2014}}%
}]{%
mohamed_improved_2014}
\APACinsertmetastar {%
mohamed_improved_2014}%
\begin{APACrefauthors}%
Mohamed, H\BPBI A\BPBI H.%
\end{APACrefauthors}%
\unskip\
\newblock
\APACrefYear{2014}.
\unskip\
\newblock
\APACrefbtitle {Improved robust adaptive control of high-order nonlinear
  systems with guaranteed performance} {Improved robust adaptive control of
  high-order nonlinear systems with guaranteed performance}\
  \APACtypeAddressSchool {M.{Sc}}{}{}.
\unskip\
\newblock
\APACaddressSchool {}{King Fahd University Of Petroleum \& Minerals}.
\PrintBackRefs{\CurrentBib}

\bibitem [\protect \citeauthoryear {%
Olfati-Saber%
, Fax%
\BCBL {}\ \BBA {} Murray%
}{%
Olfati-Saber%
\ \protect \BOthers {.}}{%
{\protect \APACyear {2007}}%
}]{%
olfati-saber_consensus_2007}
\APACinsertmetastar {%
olfati-saber_consensus_2007}%
\begin{APACrefauthors}%
Olfati-Saber, R.%
, Fax, J\BPBI A.%
\BCBL {}\ \BBA {} Murray, R\BPBI M.%
\end{APACrefauthors}%
\unskip\
\newblock
\APACrefYearMonthDay{2007}{}{}.
\newblock
{\BBOQ}\APACrefatitle {Consensus and cooperation in networked multi-agent
  systems} {Consensus and cooperation in networked multi-agent systems}.{\BBCQ}
\newblock
\APACjournalVolNumPages{Proceedings of the IEEE}{95}{1}{215--233}.
\PrintBackRefs{\CurrentBib}

\bibitem [\protect \citeauthoryear {%
Olfati-Saber%
\ \BBA {} Murray%
}{%
Olfati-Saber%
\ \BBA {} Murray%
}{%
{\protect \APACyear {2004}}%
}]{%
olfati-saber_consensus_2004}
\APACinsertmetastar {%
olfati-saber_consensus_2004}%
\begin{APACrefauthors}%
Olfati-Saber, R.%
\BCBT {}\ \BBA {} Murray, R\BPBI M.%
\end{APACrefauthors}%
\unskip\
\newblock
\APACrefYearMonthDay{2004}{}{}.
\newblock
{\BBOQ}\APACrefatitle {Consensus problems in networks of agents with switching
  topology and time-delays} {Consensus problems in networks of agents with
  switching topology and time-delays}.{\BBCQ}
\newblock
\APACjournalVolNumPages{IEEE Transactions on Automatic
  Control}{49}{9}{1520--1533}.
\PrintBackRefs{\CurrentBib}

\bibitem [\protect \citeauthoryear {%
Poggio%
\ \BBA {} Girosi%
}{%
Poggio%
\ \BBA {} Girosi%
}{%
{\protect \APACyear {1990}}%
}]{%
poggio_regularization_1990}
\APACinsertmetastar {%
poggio_regularization_1990}%
\begin{APACrefauthors}%
Poggio, T.%
\BCBT {}\ \BBA {} Girosi, F.%
\end{APACrefauthors}%
\unskip\
\newblock
\APACrefYearMonthDay{1990}{}{}.
\newblock
{\BBOQ}\APACrefatitle {Regularization algorithms for learning that are
  equivalent to multilayer networks} {Regularization algorithms for learning
  that are equivalent to multilayer networks}.{\BBCQ}
\newblock
\APACjournalVolNumPages{Science}{247}{4945}{978--982}.
\PrintBackRefs{\CurrentBib}

\bibitem [\protect \citeauthoryear {%
Qu%
}{%
Qu%
}{%
{\protect \APACyear {2009}}%
}]{%
Qu2009}
\APACinsertmetastar {%
Qu2009}%
\begin{APACrefauthors}%
Qu, Z.%
\end{APACrefauthors}%
\unskip\
\newblock
\APACrefYear{2009}.
\newblock
\APACrefbtitle {Cooperative control of dynamical systems: applications to
  autonomous vehicles} {Cooperative control of dynamical systems: applications
  to autonomous vehicles}.
\newblock
\APACaddressPublisher{}{Springer Science \& Business Media}.
\PrintBackRefs{\CurrentBib}

\bibitem [\protect \citeauthoryear {%
Ren%
\ \BBA {} Beard%
}{%
Ren%
\ \BBA {} Beard%
}{%
{\protect \APACyear {2008}}%
}]{%
ren_distributed_2008}
\APACinsertmetastar {%
ren_distributed_2008}%
\begin{APACrefauthors}%
Ren, W.%
\BCBT {}\ \BBA {} Beard, R\BPBI W.%
\end{APACrefauthors}%
\unskip\
\newblock
\APACrefYear{2008}.
\newblock
\APACrefbtitle {Distributed consensus in multi-vehicle cooperative control}
  {Distributed consensus in multi-vehicle cooperative control}.
\newblock
\APACaddressPublisher{}{Springer}.
\PrintBackRefs{\CurrentBib}

\bibitem [\protect \citeauthoryear {%
Shahvali%
\ \BBA {} Askari%
}{%
Shahvali%
\ \BBA {} Askari%
}{%
{\protect \APACyear {2016}}%
}]{%
shahvali2016cooperative}
\APACinsertmetastar {%
shahvali2016cooperative}%
\begin{APACrefauthors}%
Shahvali, M.%
\BCBT {}\ \BBA {} Askari, J.%
\end{APACrefauthors}%
\unskip\
\newblock
\APACrefYearMonthDay{2016}{}{}.
\newblock
{\BBOQ}\APACrefatitle {Cooperative adaptive neural partial tracking errors
  constrained control for nonlinear multi-agent systems} {Cooperative adaptive
  neural partial tracking errors constrained control for nonlinear multi-agent
  systems}.{\BBCQ}
\newblock
\APACjournalVolNumPages{International Journal of Adaptive Control and Signal
  Processing}{}{}{}.
\PrintBackRefs{\CurrentBib}

\bibitem [\protect \citeauthoryear {%
Shen%
\ \BBA {} Shi%
}{%
Shen%
\ \BBA {} Shi%
}{%
{\protect \APACyear {2016}}%
}]{%
shen2016output}
\APACinsertmetastar {%
shen2016output}%
\begin{APACrefauthors}%
Shen, Q.%
\BCBT {}\ \BBA {} Shi, P.%
\end{APACrefauthors}%
\unskip\
\newblock
\APACrefYearMonthDay{2016}{}{}.
\newblock
{\BBOQ}\APACrefatitle {Output Consensus Control of Multiagent Systems With
  Unknown Nonlinear Dead Zone} {Output consensus control of multiagent systems
  with unknown nonlinear dead zone}.{\BBCQ}
\newblock
\APACjournalVolNumPages{IEEE Transactions on Systems, Man, and Cybernetics:
  Systems}{46}{10}{1329--1337}.
\PrintBackRefs{\CurrentBib}

\bibitem [\protect \citeauthoryear {%
Shen%
, Shi%
, Shi%
\BCBL {}\ \BBA {} Zhang%
}{%
Shen%
\ \protect \BOthers {.}}{%
{\protect \APACyear {2016}}%
}]{%
shen2016adaptive}
\APACinsertmetastar {%
shen2016adaptive}%
\begin{APACrefauthors}%
Shen, Q.%
, Shi, P.%
, Shi, Y.%
\BCBL {}\ \BBA {} Zhang, J.%
\end{APACrefauthors}%
\unskip\
\newblock
\APACrefYearMonthDay{2016}{}{}.
\newblock
{\BBOQ}\APACrefatitle {Adaptive Output Consensus with Saturation and Dead-zone
  and its Application} {Adaptive output consensus with saturation and dead-zone
  and its application}.{\BBCQ}
\newblock
\APACjournalVolNumPages{IEEE Transactions on Industrial Electronics}{}{}{}.
\PrintBackRefs{\CurrentBib}

\bibitem [\protect \citeauthoryear {%
Theodoridis%
, Boutalis%
\BCBL {}\ \BBA {} Christodoulou%
}{%
Theodoridis%
\ \protect \BOthers {.}}{%
{\protect \APACyear {2012}}%
}]{%
theodoridis_direct_2012}
\APACinsertmetastar {%
theodoridis_direct_2012}%
\begin{APACrefauthors}%
Theodoridis, D\BPBI C.%
, Boutalis, Y\BPBI S.%
\BCBL {}\ \BBA {} Christodoulou, M\BPBI A.%
\end{APACrefauthors}%
\unskip\
\newblock
\APACrefYearMonthDay{2012}{}{}.
\newblock
{\BBOQ}\APACrefatitle {Direct adaptive neuro-fuzzy trajectory tracking of
  uncertain nonlinear systems} {Direct adaptive neuro-fuzzy trajectory tracking
  of uncertain nonlinear systems}.{\BBCQ}
\newblock
\APACjournalVolNumPages{International Journal of Adaptive Control and Signal
  Processing}{26}{7}{660--688}.
\PrintBackRefs{\CurrentBib}

\bibitem [\protect \citeauthoryear {%
Yang%
, Ge%
, Wang%
, Li%
\BCBL {}\ \BBA {} Hua%
}{%
Yang%
\ \protect \BOthers {.}}{%
{\protect \APACyear {2015}}%
}]{%
yang_adaptive_2015}
\APACinsertmetastar {%
yang_adaptive_2015}%
\begin{APACrefauthors}%
Yang, Y.%
, Ge, C.%
, Wang, H.%
, Li, X.%
\BCBL {}\ \BBA {} Hua, C.%
\end{APACrefauthors}%
\unskip\
\newblock
\APACrefYearMonthDay{2015}{}{}.
\newblock
{\BBOQ}\APACrefatitle {Adaptive neural network based prescribed performance
  control for teleoperation system under input saturation} {Adaptive neural
  network based prescribed performance control for teleoperation system under
  input saturation}.{\BBCQ}
\newblock
\APACjournalVolNumPages{Journal of the Franklin Institute}{352}{5}{1850--1866}.
\PrintBackRefs{\CurrentBib}

\bibitem [\protect \citeauthoryear {%
H.~Zhang%
\ \BBA {} Lewis%
}{%
H.~Zhang%
\ \BBA {} Lewis%
}{%
{\protect \APACyear {2012}}%
}]{%
zhang_adaptive_2012}
\APACinsertmetastar {%
zhang_adaptive_2012}%
\begin{APACrefauthors}%
Zhang, H.%
\BCBT {}\ \BBA {} Lewis, F\BPBI L.%
\end{APACrefauthors}%
\unskip\
\newblock
\APACrefYearMonthDay{2012}{}{}.
\newblock
{\BBOQ}\APACrefatitle {Adaptive cooperative tracking control of higher-order
  nonlinear systems with unknown dynamics} {Adaptive cooperative tracking
  control of higher-order nonlinear systems with unknown dynamics}.{\BBCQ}
\newblock
\APACjournalVolNumPages{Automatica}{48}{7}{1432--1439}.
\PrintBackRefs{\CurrentBib}

\bibitem [\protect \citeauthoryear {%
L.~Zhang%
, Hua%
\BCBL {}\ \BBA {} Guan%
}{%
L.~Zhang%
\ \protect \BOthers {.}}{%
{\protect \APACyear {2016}}%
}]{%
zhang2016distributed}
\APACinsertmetastar {%
zhang2016distributed}%
\begin{APACrefauthors}%
Zhang, L.%
, Hua, C.%
\BCBL {}\ \BBA {} Guan, X.%
\end{APACrefauthors}%
\unskip\
\newblock
\APACrefYearMonthDay{2016}{}{}.
\newblock
{\BBOQ}\APACrefatitle {Distributed output feedback consensus tracking
  prescribed performance control for a class of non-linear multi-agent systems
  with unknown disturbances} {Distributed output feedback consensus tracking
  prescribed performance control for a class of non-linear multi-agent systems
  with unknown disturbances}.{\BBCQ}
\newblock
\APACjournalVolNumPages{IET Control Theory \& Applications}{}{}{}.
\PrintBackRefs{\CurrentBib}

\end{thebibliography}

\clearpage
\section*{AUTHOR INFORMATION}
\vspace{30pt}
{\bf Hashim A. Hashim} is a Ph.D. candidate and a Teaching and Research Assistant in Robotics and Control, Department of Electrical and Computer Engineering at the University of Western Ontario, Canada.\\
His current research interests include Stochastic and deterministic filters on SO(3) and SE(3), control of multi-agent systems, control applications and optimization techniques.\\
\underline{Contact Information}: \href{mailto:hmoham33@uwo.ca}{hmoham33@uwo.ca}.
\vspace{40pt}

{\bf Sami El Ferik} is an Associate Professor in Control and Instrumentation, Department of Systems Engineering, at KFUPM. He obtained his B.Sc. in Electrical Engineering from Laval University, Quebec, Canada, and M.Sc. and Ph.D. both in Electrical and Computer Engineering from Polytechnique Montreal, Canada. After the completion of his Ph.D. and Post-doctor positions, he worked with Pratt and Whitney Canada at the Research and Development Center of Systems, Controls, and Accessories. \\
His research interests are in sensing, monitoring, multi-agent systems and nonlinear control with strong multidisciplinary research and applications.\\
\vspace{40pt}

{\bf Frank L. Lewis} National Academy of Inventors, Fellow IEEE, Fellow IFAC, PE Texas, U.K. Chartered Engineer, is a UTA Distinguished Scholar Professor, UTA Distinguished Teaching Professor, and Moncrief-O'Donnell Chair at the University of Texas at Arlington Research Institute. He obtained the Bachelor's Degree in Physics/EE and the MSEE at Rice University, the M.S. in Aeronautical Engineering from the University of West Florida, and the Ph.D. at Ga. Tech. He is author of seven U.S. patents, numerous journal papers, and 14 books.\\

\end{document}